\documentclass[10pt,a4paper]{article}   

\usepackage[margin=1.25in]{geometry}
\usepackage{hyperref}
\usepackage[pdftex]{graphicx}
\usepackage{amsmath,amssymb,amsfonts}
\usepackage[textsize=small,textwidth=1.3in]{todonotes}
\usepackage[english]{babel}
\usepackage{enumerate}
\usepackage{multicol}
\usepackage{subcaption}
\usepackage[ruled]{algorithm2e}
\usepackage{graphicx}
\usepackage{subcaption}

\usepackage{algorithmic}

\usepackage[utf8]{inputenc}
\usepackage[margin=1.25in]{geometry}
\usepackage[T1]{fontenc}
\usepackage[english]{babel}
\usepackage{color}
\usepackage[textsize=small]{todonotes}

\usepackage{amsfonts}
\usepackage{amsmath,stackrel,bm}
\usepackage{amsthm}
\usepackage{enumerate}
\usepackage{cite}
\usepackage{multicol}
\usepackage{multirow}
\usepackage{fullpage}
\usepackage{hyperref}
\usepackage{soul}

\usepackage[all,dvips]{xy}
\usepackage{graphicx}
\usepackage{tikz}
\usetikzlibrary{matrix,arrows,decorations.pathmorphing}

\theoremstyle{plain}
\newtheorem{thm}{Theorem}[section]
\numberwithin{thm}{section}
\newtheorem{prop}[thm]{Proposition}
\newtheorem{lema}[thm]{Lemma}
\newtheorem{cor}[thm]{Corollary}
\newtheorem{defi}[thm]{Definition}
\newtheorem{remark}[thm]{Remark}
\newtheorem{ex}[thm]{Example}

\newtheorem*{keywords}{Keywords}
\newtheorem*{PACS}{Mathematics Subject Classification (2010)}

\newcommand{\mS}{{\mathcal S}}
\newcommand{\mA}{{\mathcal A}}
\newcommand{\ba}{{\bf a}}
\newcommand{\zz}{{\mathbb Z}}
\newcommand{\bz}{{\bf 0}}
\newcommand{\bb}{{\bf b}}
\newcommand{\bc}{{\bf c}}
\newcommand{\mL}{{\mathcal L}}
\newcommand{\eqdef}{\stackrel{\text{def}}{=}}
\newcommand{\n}{{\mathbb N}}
\newcommand{\balpha}{\bm{\alpha}}
\newcommand{\K}{{\mathbb K}}
\newcommand{\bbeta}{\bm{\beta}}

\newcommand{\eproof}{\hfill$\square$}
\newcommand{\C}{{\cal C}}
\newcommand{\mT}{{\mathcal T}}
\newcommand{\blambda}{\bm{\lambda}}
\newcommand{\bnu}{\bm{\nu}}
\newcommand{\bgamma}{\bm{\gamma}}
\newcommand{\bdelta}{\bm{\delta}}
\newcommand{\bx}{{\bf x}}

\usepackage{authblk}

\begin{document}

\title{Factorizations of the same length in abelian  monoids\thanks{This work was partially supported by the Spanish MICINN PID2019-105896GB-I00 and MASCA (ULL Research Project)}}
\author[1]{Evelia R. Garc\'ia Barroso }
\author[1]{Ignacio Garc\'ia-Marco}
\author[1]{Irene M\'arquez-Corbella}
\affil[1]{
Departamento de Matem\'aticas, Estad\'istica e I.O., Universidad de La Laguna, Apdo. Correos 460, 38200, La Laguna, Tenerife, Spain.

Email: 
 \href{mailto:ergarcia@ull.es}{ergarcia@ull.es},
 \href{mailto:iggarcia@ull.es}{iggarcia@ull.es} and 
\href{mailto:imarquec@ull.es}{imarquec@ull.es} 
}
%
%

%
%

\maketitle

\begin{abstract}
Let $\mS \subseteq \zz^m \oplus T$ be a finitely generated and reduced  monoid. In this paper we develop a general strategy to study the set of elements in $\mathcal S$ having at least two factorizations of the same length, namely the ideal $\mathcal L_{\mathcal S}$.
To this end, we work with a certain (lattice) ideal associated to the monoid $\mathcal S$.
Our study can be seen as a new approach generalizing  \cite{chapman:2011}, which only studies the case of numerical semigroups. When $\mS$ is a numerical semigroup we give three main results: (1) we compute explicitly a set of generators of the ideal $\mathcal L_{\mathcal S}$ when $\mathcal S$ is minimally generated by an almost arithmetic sequence;  (2) we provide an infinite family of numerical semigroups such that 
$\mathcal L_{\mathcal S}$ is a principal ideal; (3) we classify the computational problem of determining the largest integer not in $\mathcal L_{\mathcal S}$ as an $\mathcal{NP}$-hard problem.

\end{abstract}


\section{Introduction}
\label{intro}

Let $\mS$ be an abelian monoid, it is well-known that $\mS$ can be embedded in a group $\mathcal G$ if and only if $\mS$ is cancellative (see \cite{clifford:1967}). The usual procedure for doing this is by considering the so called {\it group of quotients} of $\mS$. That is, the abelian  group $\mathcal G = (\mS \times \mS)/ \sim$, where $(a,b) \sim (c,d)$ if and only if $a + d = b + c$.  This group $\mathcal G$ contains $\mS$ via the embedding $s \mapsto (s,\mathbf 0)$. 

An abelian monoid $\mS$ is {\it finitely generated} if there exist some $\ba_{1},\ldots,\ba_{n}\in \mS$ such that $$\mS =\left\{ \lambda_1 \ba_1 + \cdots + \lambda_n \ba_n \mid \lambda_1, \ldots, \lambda_n \in \n\right\},$$ in which case we will put $\mS = \left\langle \ba_1, \ldots, \ba_n\right\rangle$. One has that $\mathcal G$ is finitely generated whenever $\mS$ so is.

Therefore, if $\mS$ is an abelian, cancellative, finitely generated monoid then
$$\mS = \left\langle \ba_1, \ldots, \ba_n\right\rangle
\subseteq \mathcal G \simeq \mathbb Z^m \oplus T,$$
where  $T$ is a finite abelian group
and $m=\mathrm{rank}(\mathcal G)$ is the rank of $\mathcal G$.
If $\mathcal G$ is torsion-free, then $T=\{\bz\}$ and $\mS$ is called an {\em affine monoid}.

The monoid $\mS\subseteq \mathbb Z^m \oplus T$ is {\em reduced} if the only invertible element is the neutral element of $\mS$ or, equivalently, if $\mS \cap (-\mS) = \{\bf{0}\}$. If $\mS$ is reduced, then $\mS$ has a unique minimal (with respect to the inclusion) set of generators, which coincides with the {\em set of atoms} or {\em irreducible elements} of $\mS$. We will refer to this set of generators as  {\em the minimal set of generators} of the monoid $\mS$.

Unless otherwise stated, when we write $\mS = \left\langle \ba_1, \ldots, \ba_n\right\rangle\subseteq \mathbb Z^m \oplus T$ for a reduced monoid, we are assuming that $\mS$ is an abelian, cancellative monoid that can be embedded in $\mathbb Z^m \oplus T$, being $T$ a finite abelian group, 
and $\mathcal A = \left\{ \ba_1, \ldots, \ba_n \right\}$ is the minimal set of generators of $\mS$. These monoids provide a powerful interface between Combinatorics and Algebraic Geometry since they constitute a combinatorial tool for studying lattice ideals and Toric Geometry (see, e.g. \cite{sturmfels:1996,BGT:2002,cox:2011,villarreal:2015}).

 Now, consider a reduced monoid $\mS = \langle \ba_1,\ldots,\ba_n \rangle \subseteq \zz^m\oplus T$. For any  $\bb \in \mS$ there exists an $n$-tuple $\bm{\lambda} = (\lambda_1, \ldots, \lambda_n)\in \n^n$ such that 
$\bb = \lambda_1 \ba_1 + \cdots + \lambda_n \ba_n$.
In this case we say that $\bm{\lambda} = (\lambda_1, \ldots, \lambda_n)$ is a {\it factorization} of $\bb$ in $\mS$ of {\it length} $\ell(\bm{\lambda}) := \lambda_1 + \cdots + \lambda_n$. 
Now we define the set $\mL_{\mS}$ of elements in $\mS$ having (at least) two factorizations of the same length, i.e.,
\[ \mL_{\mS} \eqdef \left\{ \bb \in \mS \mid \bb \hbox{ has two different factorizations of the same length}\right\}.\]

In this paper we investigate the set $\mL_{\mS}$. In the particular setting that $\mS$ is a numerical semigroup this problem was addressed in \cite{chapman:2011}. Numerical semigroups provide an interesting family of reduced monoids with $T=\{0\}$ (affine monoids). More precisely, a {\it numerical semigroup} is a submonoid of $\n$ with finite complement over $\n$  (for a thorough study of numerical semigroups we refer the reader to \cite{Assi:2016,alfonsin:1996}). In \cite{chapman:2011}, the authors prove that given a numerical semigroup $\mS = \langle a_1,\ldots,a_n\rangle \subseteq \n$, then $\mL_{\mS} = \emptyset$ if and only if $n = 2$, and describe $\mL_{\mS}$ when $n = 3$.

This paper goes further into the study of factorization properties of reduced monoids by means of their corresponding lattice ideal.  See \cite{geroldinger:2006} for a general reference in the theory of non-unique factorization domains and monoids. For a recent account of the progress of factorization invariants in affine monoids, we refer the reader to the recent papers \cite{GOW:2019,GZ:2020} and the references therein.

\subsection*{\bf Outline of the article}
Section \S\ref{Sec:Apery} is devoted to the study of the Ap\'ery set of a reduced monoid with respect to a finite set $B = \{\bb_1,\ldots,\bb_s\} \subseteq \mS$, that is, the set
\[  \mathrm{Ap}_{\mS}(B) = \left\{ x\in \mS \mid x-\bb_i \notin \mS, \; 1\leq i \leq s \right\}. \]
Although we use Ap\'ery sets in the other sections, we believe that the results in this section are  interesting in their own. In Theorem \ref{thm:apery}, we present how to compute an Ap\'ery set by means of the (lattice) ideal $I_\mS$ of the monoid and a factorization of the elements of $B$. This result provides an alternative to \cite[Theorem 8]{marquez-campos:2015}. Then, in Theorem \ref{teo:apfinito}, we characterize when this Ap\'ery set is finite, it turns out that $\mathrm{Ap}_{\mS}(B)$ is finite if and only if the union of $\{\mathbf{0}\}$ and the  ideal generated by $B$ form a reduced monoid. We prove that this is also equivalent to the fact that the cone defined by $B$  (see Definition \ref{def:polyhedral}) coincides with the one defined by $\mS$.

The main results of the paper are in Section \S \ref{sec:samelenght}, where we develop a general strategy to study $\mL_\mS$. In Proposition \ref{prop:dosescrituras} we describe how to obtain a finite set of generators of the ideal $\mL_\mS$ by means of the lattice ideal $I_{\tilde{\mS}}$ of the monoid 
\begin{equation} 
\label{eq:Stilde}
\tilde{\mS}= \left\langle (\ba_1, 1) , (\ba_2, 1), \ldots , (\ba_n,1)\right\rangle \subseteq \zz^{m+1} \oplus T, 
\end{equation}
associated with $\mS$.

As a consequence, in Theorem \ref{thm:main}, we describe $\mS \setminus \mL_\mS$ as an Ap\'ery set and, thus, the techniques developed in \S\ref{Sec:Apery} apply here. In particular, using Theorem \ref{teo:apfinito}, we describe when $\mL_\mS \cup \{\bz\}$ is a reduced monoid or, equivalently, when $\mS \setminus \mL_\mS$ is a finite set (see Corollary \ref{cor:maincoro}).
 In the last part of this section we apply our results to the particular context of numerical semigroups and provide alternative  proofs of the results of \cite{chapman:2011} mentioned above. 

In Section \S \ref{sec:catenaryDegree} we study the notion of equal catenary degree of a reduced monoid. Equal catenary degrees have been studied since $2006$, see for example \cite{foroutan:2006,geroldinger:2006,hassler:2009,blanco:2011,philipp:2015,geroldinger:2019} and the references therein. Our main result in this section is Theorem \ref{th:catenary}, where we prove that $\rm c_{\rm eq}(\mS)$ equals the maximum degree of a minimal generator of $I_{\tilde{\mS}}$. In particular we improve \cite[Proposition 4.4.3]{blanco:2011} and recover \cite[Lemma 6]{gonzalez2019monotone}. Then, applying to $I_{\tilde{\mS}}$ the upper bound for the Castelnuovo-Mumford regularity of projective monomial curves given by L'vovsky in \cite{Lvovsky}, we obtain in Theorem \ref{th:cotacurvas} an upper bound for the equal catenary degree of any numerical semigroup.

Section \S \ref{sec:arithmetic} is devoted to prove Theorem \ref{prop:almostarithmetic-A}, where we provide an explicit set of generators of the ideal $\mL_\mS$ when $\mS$ is minimally generated by an almost arithmetic sequence. By almost arithmetic sequence we mean a set $\{m_1,\ldots,m_n,b\}$, where $m_1<\ldots<m_n$ is an arithmetic sequence of positive integers and $b$ is any positive integer. The key idea to prove these results is to use \cite[Theorem 2.2]{bermejo:2017}. There, the authors describe a set of generators of the ideal of some projective monomial curves, which, in this context coincide with the toric ideal of $\tilde{\mS}$, and then we apply Proposition \ref{prop-main-dS} with this set of generators.

In Section \S \ref{sec:principal} we address the question of characterizing when $\mL_\mS$ is a principal ideal. We give a partial answer to this question by providing in Corollary \ref{cor:uniqueBetti} an infinite family of numerical semigroups such that $\mL_\mS$ is a principal ideal. This family consists of shiftings of numerical semigroups with a unique Betti element (a family of semigroups studied in \cite{garcia-sanchez:2012}), and generalizes three generated numerical semigroups. As an intermediate result, in Proposition \ref{prop:descuniquebetti} we describe an explicit set of generators of the toric ideal of a family of numerical semigroups which turn to be semigroups with a single Betti minimal element (a family of semigroups studied in \cite{garcia-sanchez:2019}).

When $\mS \subseteq \n$ is a numerical semigroup and $\mL_\mS$ is not empty, then $\n ~\setminus~\mL_\mS$ is a finite set. In Section \S \ref{sec:NP} we classify the computational problem of determining the largest integer not in $\mL_\mS$ as an $\mathcal{NP}$-hard problem. We derive this result by restating the proof of the $\mathcal{NP}$-hardness of the Frobenius problem in \cite{alfonsin:1996}  and  some (easy) considerations. The same ideas also allow us to derive that, for a bounded value $k \in \zz^+$, computing the largest element in a numerical semigroup with at least $k$ different factorizations (or at least $k$ different factorizations of the same length) is $\mathcal{NP}$-hard.

\section{Ap\'ery sets of reduced monoids}
\label{Sec:Apery}

Let $\mS = \left\langle \ba_1, \ldots , \ba_n \right\rangle\subseteq \zz^m \oplus T$ be a reduced monoid and consider a finite set of nonzero elements $B=\left\{ \bb_1, \ldots, \bb_s \right\} \subseteq \mS \setminus \{\bz\}$.  We define the Ap\'ery set of $\mS$ with respect to $B$ as 
\begin{equation}\label{def:apery} \mathrm{Ap}_{\mS}(B) = \left\{ x\in \mS \mid x-\bb_i \notin \mS, \; 1\leq  i \leq s \right\}. \end{equation}

In this section we study ${\rm Ap}_{\mS}(B)$ and provide Theorem \ref{thm:apery} and Theorem \ref{teo:apfinito} as the main results. In the first we describe Ap\'ery sets in terms of the degrees of the elements of a certain basis of a $\mathbb K$-vector space. In the second one we characterize when Ap\'ery sets are finite.

The problem of computing the Ap\'ery set of an affine monoids has been studied in \cite{pison:2003,marquez-campos:2015}. In \cite{marquez-campos:2015} the authors provide a method to compute the Ap\'ery set of 
an affine semigroup based on Gr\"obner basis computations. Our Theorem \ref{thm:apery} is more general, since we do not require that $\mathcal G$ is torsion-free, but it is inspired by \cite[Theorem 8]{marquez-campos:2015}. However, even in the affine monoid setting, the main differences are: 
In \cite{marquez-campos:2015}, the authors require an extra hypothesis implying that the Ap\'ery set is finite that we do not assume. In Proposition \ref{teo:apfinito} we prove that this extra hypothesis characterizes when ${\rm Ap}_{\mS}(B)$ is finite.  Another difference is that our result does not need any choice of a monomial order.
A third difference is that Theorem \ref{thm:apery} requires a factorization of $\bb_1,\ldots,\bb_r$, while in \cite{marquez-campos:2015} they do not require so. This is 
not a big limitation for us, since we are applying this result in Section~\ref{sec:samelenght} in a context where we already know a factorization of the elements of $B$.

To state and prove Theorem \ref{thm:apery}, first we will introduce some basic notions on lattice ideals.
Let $\K$ be a field, we denote by $\K[\mathbf{x}] = \K[x_1, \ldots, x_n]$ the ring of polynomials in the variables $x_1,\ldots,x_n$ 
with coefficients in $\K$.  We write a monomial in $\K[\mathbf{x}]$ as 
$$\mathbf{x}^{\balpha} = x_1^{\alpha_1} \cdots x_n^{\alpha_n} \hbox{ with } \balpha = (\alpha_1, \ldots, \alpha_n) \in \n^n.$$

A reduced monoid $\mS = \left\langle \ba_1, \ldots, \ba_n\right\rangle  \subseteq \zz^m \oplus T$ induces a grading in $\mathbb \K[\mathbf{x}]$ given by 
\[\deg_\mS(\mathbf{x}^{\balpha}) = \sum_{i=1}^n \alpha_i \ba_i, \hbox{ for } \balpha = (\alpha_1, \ldots, \alpha_n) \in \n^n, \]
and called $\mS$-degree.

A polynomial $f\in \mathbb \K[\mathbf{x}]$ is $\mS$-homogeneous if all its monomials have the same $\mS$-degree. Moreover, an ideal is $\mS$-homogeneous if it is generated by $\mS$-homogeneous polynomials.

Associated to $\mS$, we have the monoid algebra $\mathbb K[\mS] = \mathbb K[\mathbf{t}^s \mid s\in \mS]$. Consider the epimorphism of $\mathbb K$-algebras:
\begin{equation}
\label{toric:morphism}
\begin{array}{cccl}
\varphi: & \K[\mathbf{x}] & \longrightarrow & \K[\mS],\\
& x_i & \longmapsto & \mathbf{t}^{\ba_i}  
\end{array}
\end{equation}
the \emph{lattice ideal} of $\mS$ is $I_\mS = \ker(\varphi)$.

It turns out that $\mathbb K[\mS]$ is an integral domain if and only if the group of quotients of $\mS$ is torsion-free or, equivalently, if $\mS$ is an affine monoid. In this case $\mathbb K[\mS]$ becomes a subalgebra of the Laurent polynomial ring $\mathbb K[t_1^{\pm1},\ldots, t_m^{\pm1}]$. On the other hand, lattice prime ideals are called toric ideals. 
Hence the ideal $I_{\mS}$ is toric if and only if $K[\mS]$ is an integral domain and, thus, this is equivalent to $\mS$ is an affine monoid. 

\begin{remark}
\label{rem:binomios}
The lattice ideal $I_\mS$ has been thoroughly studied in the literature 
(see, e.g., \cite{sturmfels:1996,villarreal:2015}). For example, it is well known that $I_\mS$ is an $\mS$-homogeneous binomial ideal (it is generated by differences of monomials). We have that $\mathbf x^{\balpha}- \mathbf x^{\bbeta}\in I_{\mS}$ if and only if $\deg_{\mS}(\mathbf x^{\balpha}) = \deg_{\mS}(\mathbf x^{\bbeta})$; as consequence
\begin{equation}
\label{eq:binomios}
I_{\mS} = \left\langle \mathbf x^{\balpha} - \mathbf x^{\bbeta} \mid \deg_{\mS}(\mathbf x^{\balpha}) = \deg_{\mS}(\mathbf x^{\bbeta})\right\rangle. \end{equation}

 Moreover, $I_\mS$ is of height ${\rm ht}(I_\mS) = n - {\rm rank}(\mathcal G)$, where $\mathcal G$ is the group of quotients of $\mS$. Equivalently ${\rm rank}(\mathcal G)={\rm rank}(A)$, where $A$
is the $m \times n$ matrix with columns $\pi(\ba_1),\ldots,\pi(\ba_n) \in \mathbb Z^m$, being
$\pi$  the canonical projection 
\begin{equation}
\label{canonicalProj}
\begin{array}{cccc}\pi : & \mathbb Z^m \oplus T & \longrightarrow & \mathbb Z^m\\ 
& (x,t) & \longmapsto & x\end{array}.
\end{equation}

Consider the group homomorphism
$\rho:  \mathbb Z^n  \longrightarrow  \mathbb Z^m$ such that $\rho(\mathbf e_i) = \ba_i,$
where $\{\mathbf e_1,\ldots,\mathbf e_n\}$ is the canonical basis of $\mathbb Z^n$. From (\ref{eq:binomios}) one deduces that $I_{\mS} =  \left\langle \mathbf x^{\balpha} - \mathbf x^{\bbeta} \mid {\balpha} - {\bbeta} \in \ker(\rho) \right\rangle$. Hence, this ideal can be computed in the following way:
Compute a generating set of the kernel of $\rho$, i.e. $\ker(\rho) = \left\langle \bgamma_1, \ldots, \bgamma_{t}\right\rangle\subseteq \mathbb Z^n$ and write every element $\bgamma_i\in \mathbb Z^n$ as $\bgamma_i = \bgamma_i^+ - \bgamma_i^-$ with $\bgamma_i^+, \bgamma_i^- \in \mathbb N^n$. Then,
\begin{equation} \label{eq:saturation} I_{\mS} = \left(\left\langle \mathbf x^{\bgamma_i^+}- \mathbf x^{\bgamma_i^+} \mid 1 \leq i \leq t\right\rangle : \left(x_1 \cdots x_n\right)^{\infty}\right). \end{equation}
Recall that, if $J\subseteq \mathbb K[\mathbf x]$ is an ideal then \[J:f^{\infty} = \left\{ g\in \mathbb K[\mathbf x] \mid  \hbox{ there is } k \geq 1 \hbox{ such that } gf^k \in J\right\}\] is again an ideal of $\mathbb K[\mathbf x]$.

The expression (\ref{eq:saturation}) provides a method for computing a set of generators of $I_\mS$; for improvements of this method see, e.g., \cite{BLR,CTV2017,Hemmecke2009}.

Moreover, since $\mS$ is reduced, a graded version of Nakayama's lemma holds. As a consequence, all minimal sets of binomial generators  of $I_{\mS}$ have the same number of elements and the same $\mS$-degrees.
\end{remark}

 Consider now $B = \{\bb_1,\ldots,\bb_s\} \subseteq \mS \setminus\{\bz\}$. Since $\bb_i \in \mS$, one can express $\bb_i =\sum_{j=1}^n \beta_{ij}\ba_j$, where $\bbeta_i = (\beta_{i1},\ldots,\beta_{in})\in \n^n$. Let $\mathbf{x}^{\bbeta_i} = x_1^{\beta_{i1}} \cdots x_n^{\beta_{in}}$ for all $i\in \{1, \ldots, s\}$.

\begin{thm}
\label{thm:apery}Let $\mS = \langle \ba_1,\ldots,\ba_n \rangle\subseteq \mathbb Z^m \oplus T$ be a reduced monoid and let $B = \{\bb_1,\ldots,\bb_s\} \subseteq \mS \setminus \{\bz\}.$ Set the monomial $\mathbf{x}^{\bbeta_i} = x_1^{\beta_{i1}} \cdots x_n^{\beta_{in}} \in \mathbb \K[\mathbf{x}],$ where $\bbeta_i = (\beta_{i1},\ldots,\beta_{in})\in \n^n$ is a factorization of $\bb_i$ for all $i\in \{1, \ldots, s\}$. If we take  a monomial $\mathbb K$-basis $D$ of $\mathbb \K[\mathbf{x}] / (I_\mS + \left\langle \mathbf{x}^{\bbeta_1}, \ldots, \mathbf{x}^{\bbeta_s}\right\rangle)$, then the mapping
\[
\begin{array}{cccc}
h: & D & \longrightarrow & \mathrm{Ap}_S(B) \\
& \mathbf{x}^{\balpha} & \longmapsto & \deg_{\mS}(\mathbf{x}^{\balpha}) = \alpha_1 \ba_1 + \cdots + \alpha_n \ba_n
\end{array}
\] 
is bijective.
\end{thm}

\proof
We start with the epimorphism presented in Equation \eqref{toric:morphism}. We observe that $\varphi$ is graded with respect to the grading ${\rm deg}_\mS(x_i) = \ba_i$ and $\deg(t^{\bb}) = \bb \in \mS$. We have that $\mathbb \K[\mathbf{x}] / I_\mS \simeq \mathbb K[\mS]$ and we denote by $\tilde{\varphi}$ the corresponding graded isomorphism of $\mathbb K$-algebras.

Now we consider  the ideal $\langle \mathbf{t}^{\bb_1}, \ldots, \mathbf{t}^{\bb_s}\rangle \cdot \K[\mS]$ generated by  $\mathbf{t}^{\bb_1}, \ldots, \mathbf{t}^{\bb_s}$  in $\K[\mS]$, and the canonical epimorphism:
\[
\begin{array}{cccc} 
e : & \mathbb K[\mS] & \longrightarrow & \mathbb K[\mS] / \langle \mathbf{t}^{\bb_1}, \ldots, \mathbf{t}^{\bb_s}\rangle  \cdot \K[\mS]\\
& \mathbf{t}^{\balpha} & \longmapsto & [\mathbf{t}^{\balpha}].
\end{array}
\]

Since $\varphi(\mathbf{x}^{\bbeta_i})= \mathbf{t}^{\bb_i}$, we have that  $\ker(e \circ \tilde{\varphi}) = (I_\mS + \langle \mathbf{x}^{\bbeta_1}, \ldots, \mathbf{x}^{\bbeta_s}\rangle) / I_\mS$. Thus, by the third isomorphism theorem, there is a graded isomorphism of $\mathbb K$-algebras
\[
\Psi: \mathbb \K[\mathbf{x}] / (I_\mS + \langle \mathbf{x}^{\bbeta_1}, \ldots, \mathbf{x}^{\bbeta_s}\rangle) \longrightarrow \mathbb \K[\mS]/ \langle \mathbf{t}^{\bb_1}, \ldots, \mathbf{t}^{\bb_s} \rangle \cdot \K[\mS].
\]

Moreover,  $\K[\mS]/ \langle \mathbf{t}^{\bb_1}, \ldots, \mathbf{t}^{\bb_s} \rangle \cdot \K[\mS]$ has a unique monomial basis, which is $\{\mathbf{t}^{\bb} \, \vert \, \bb \in {\rm Ap}_\mS(B)\}$. Finally, we observe that the image of a monomial by $\Psi$ is a monomial and hence, the image of any monomial basis $D$ of 
 $\mathbb \K[\mathbf{x}] / (I_\mS + \left\langle \mathbf{x}^{\bbeta_1}, \ldots, \mathbf{x}^{\bbeta_s}\right\rangle)$ has to be  $\{\mathbf{t}^{\bb} \, \vert \, \bb \in {\rm Ap}_\mS(B)\}$. The result follows from the fact that $\Psi$ is graded and $\Psi(\mathbf{x}^{\balpha}) = \mathbf{t}^{\deg_\mS(\mathbf{x}^{\balpha})}$.
\eproof

\medskip

Set $J := I_\mS + \left\langle \mathbf{x}^{\bbeta_1}, \ldots, \mathbf{x}^{\bbeta_s}\right\rangle$. To compute  a monomial $\mathbb K$-basis $D$ of $\mathbb \K[\mathbf{x}]/J$, it suffices to choose any monomial ordering $\succ$ in $\mathbb \K[\mathbf{x}]$, and define $D$ as the set of all the monomials not belonging to $\mathrm{in}_{\succ} (J)$, the initial ideal of $J$ with respect to $\succ$. That is,   
$$D=\left\{ \mathbf{x}^{\balpha} \mid \mathbf{x}^{\balpha} \notin \mathrm{in}_{\succ}(J)\right\}.$$
Notice that different monomial orders yield different $\mathbb K$-bases. Nevertheless, Theorem \ref{thm:apery} holds for any of these (and for any other monomial $\mathbb K$-basis).

Let us illustrate the previous result with some examples.
\begin{ex}\label{ex:apery}Let $\mS = \langle \ba_1,\ldots,\ba_5\rangle \subseteq \zz^2$ with $\ba_1 = (0,2), \ba_2 = (1,2), \ba_3 = (1,1), \ba_4 = (3,2), \ba_5 = (4,2)$ and consider the set $B = \{\bb_1, \bb_2, \bb_3\} \subseteq \mS$, where $\bb_1 = (3,6), \bb_2 = (4,4), \bb_3 = (9,6).$ A computation with any software for polynomial computations (e.g., {\sc Singular} \cite{DGPS}, CoCoA \cite{CoCoA} or {\tt Macaulay2} \cite{M2}) shows that $I_{\mS} = \langle f_1,\ldots,f_6\rangle$ with 
\[ \begin{array}{lll} f_1 = x_4^2-x_3^2x_5, & f_2 =x_3^2 x_4-x_2 x_5, & f_3 = x_2 x_4 -x_1 x_5, \\ f_4 = x_3^4-x_1 x_5, & f_5 = x_2 x_3^2-x_1 x_4, & f_6 = x_2^2-x_1 x_3^2. \end{array} \]
Let us  compute  a factorization $\bbeta_i$ of $\bb_i$ for $i  \in \{ 1,2,3\}$: 

\[ \begin{array}{lll} \bb_1 = 3 \ba_2, &  \bb_2 = \ba_2 + \ba_4, & \bb_3 = 3 \ba_4,  \end{array} \]
and set \[ \begin{array}{lll} \mathbf{x}^{\bbeta_1} = x_2^3, &  \mathbf{x}^{\bbeta_2} = x_2 x_4, & \mathbf{x}^{\bbeta_3} = x_4^3.  \end{array} \] If one 
considers $L = {\rm in}_\succ(I_\mS + \langle  x_2^3, x_2 x_4, x_4^3 \rangle),$ where $\succ$ is the weighted degree reverse lexicographic order with weights $(2,2,1,2,2)$, 
then one gets \[ L = \langle x_1^2 x_4,\, x_1 x_5,\, x_2^2,\, x_2 x_3^2,\, x_2 x_4,\, x_2 x_5^2,\, x_3^4,\, x_3^2 x_4,\, x_4^2 \rangle.\] 
Hence, the monomials which are not in $L$ form the following  monomial $\K$-basis of $\K[x_1,\ldots,x_5] / (I_\mS + \langle  x_2^3, x_2 x_4, x_4^3 \rangle)$:
\[ \begin{array}{llll} D  =  & \{x_1^{a} x_3^{c}, \ x_3^{c} x_5^a \,\vert \, a \in \n,\, c \in \{0,1,2,3\} \} & \cup  \\   
& \{x_1^{a} x_2 x_3^{c}, \ x_3^{c} x_4 x_5^a \,\vert \, a \in \n,\, c \in \{0,1\} \} & \cup  \\ 
    &   \{x_1x_4, x_2x_5, x_1 x_3x_4,x_2x_3,x_5\}. \end{array} \]
Thus, by Theorem \ref{thm:apery}, the Ap\'ery set with respect to $B$ is the infinite set
\[ \begin{array}{llll}  {\rm Ap}_\mS(B)  = 	& \left\{(i,i+2\lambda),\, (i + 4\lambda,i + 2 \lambda)\, \vert \, \lambda \in \n,\ i \in \{0,1,2,3\} \right\} \  \cup \\
											& \left\{x + \lambda (0,2)\, \vert \, \lambda \in \n,\ x \in \{(1,2), (2,3)\} \right\} \  \cup \\
											& \left\{ x + \lambda(4,2) \, \vert \, \lambda \in \n,\ x \in \{(3,2) , (4,3) \} \right\} \  \cup \\
											& \{(3,4), (5,4), (4,5), (6,5) \}. \end{array} \]
See Figure \ref{fig:apery} for a graphical representation of  ${\rm Ap}_\mS(B)$. 

\begin{figure}
\includegraphics[scale=1]{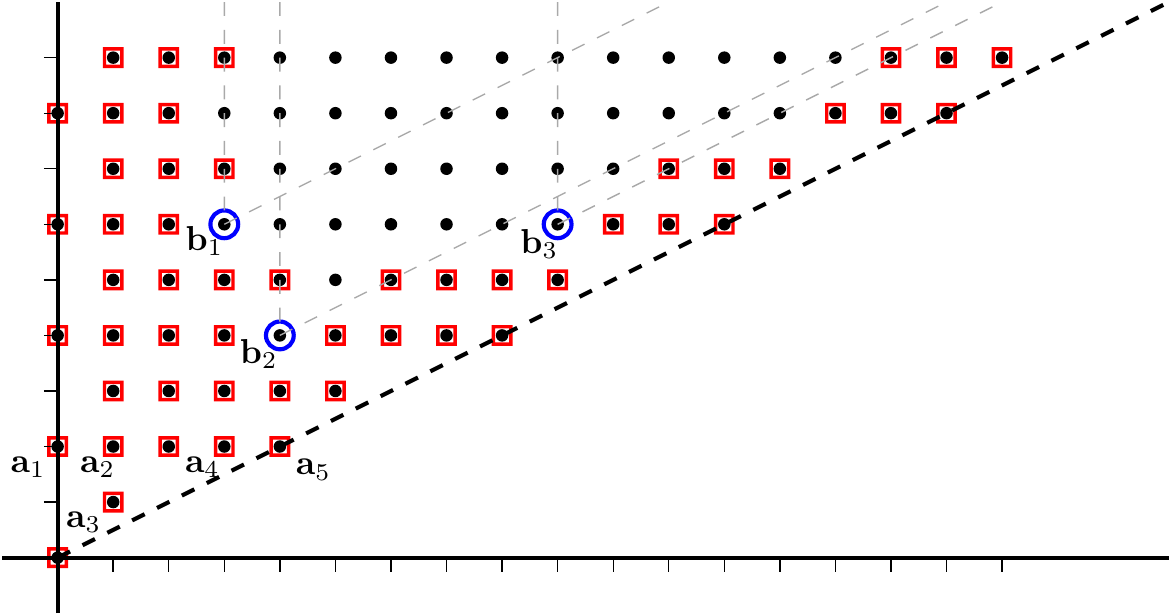} 
\caption{Ap\'ery set ${\rm Ap}_\mS(B)$ in Example \ref{ex:apery}. The dots correspond to the elements in $\mS$, the circles to the elements in $B$  and the
squares to the e\-le\-ments in ${\rm Ap}_\mS(B)$.} \label{fig:apery} 
\end{figure}
\end{ex}

\begin{ex}
\label{ex:apery2}
Let $\mS = \left\langle \ba_1, \ba_2, \ba_3 \right\rangle \subseteq \mathbb Z \oplus \mathbb Z_2$ with $\ba_1 = (2, \overline{0})$, $\ba_2 = (3, \overline{1})$ and $\ba_3 = (4, \overline{1})$ and consider the set $B=\{ \bb_1\}\subseteq \mS$ with $\bb_1 =(12, \overline{0}) $. 
One has that
$\ker(\rho) = \left\langle (1,2,-2), (0,8,-6)\right\rangle$. Thus, by Remark \ref{rem:binomios}, 
\[I_{\mS} = \langle x_1x_2^2 - x_3^2, x_2^8-x_3^6 \rangle : (x_1x_2x_3)^{\infty} = \langle x_1x_2^2-x_3^2, x_1^3-x_2^2, x_2^4-x_1^2x_3^2\rangle.\]
Let us compute a factorization $\beta_1$ of $\bb_1$, that is, $\bb_1 = 4 \ba_2$ and set $\mathbf x^{\beta_1} = x_2^4$. If one considers $L = {\rm in}_\succ(I_\mS + \langle  x_2^4 \rangle)$ where $\succ$ is the degree reverse lexicographic order, then one gets
\[L=\left\langle x_1x_2^2, x_1^3, x_3^4, x_2x_3^2, x_1^2x_2^2, x_2^4 \right\rangle.\]
Hence, the monomials which are not in $L$ form the following monomial $\mathbb K$-basis of $\K[x_1,x_2,x_3] / (I_\mS + \langle  x_2^4 \rangle)$:
\[\begin{array}{ccc}
D & = & \left\{\begin{array}{c} 
1, x_1, x_1^2, x_2, x_1x_2, x_1^2x_2, x_2^2, x_2^3, x_3, x_1x_3, x_1^2x_3, x_2x_3,x_1x_2x_3, x_1^2x_2x_3, \\
x_2^2x_3, x_2^3x_3,x_3^2, x_1x_3^2, x_2x_3^2, x_1x_2x_3^2, x_1x_2x_3^2, x_3^3, x_1x_3^3, x_2x_3^3, x_1x_2x_3^3
\end{array}\right\}\\
\end{array}.
\]
Thus, by Theorem \ref{thm:apery}, the Ap\'ery set with respect to $B$ is the finite set
\[{\rm Ap}_\mS(B)  =  \{ {\rm deg}_{\mS}(\bx^{\alpha}) \, \vert \, \bx^{\alpha} \in D\},\]
which is \[
 \{ (x, \overline{0}) \mid  x \in \{0,2,4,6,7,8,9,10,11,13,15,17\}\} \cup 
 \{ (x, \overline{1}) \mid 3 \leq x \leq 14 \}. \]
 See Figure \ref{fig:apery2} for a graphical representation of  ${\rm Ap}_\mS(B)$. 

\begin{figure}
\includegraphics[scale=1]{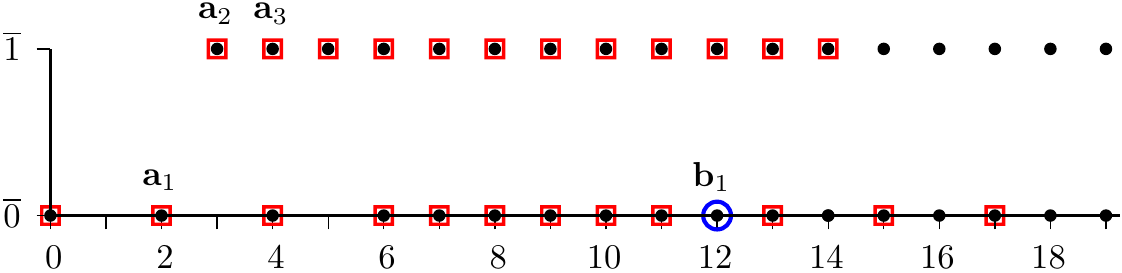} 
\caption{Ap\'ery set ${\rm Ap}_\mS(B)$ in Example \ref{ex:apery2}. The dots correspond to the elements in $\mS$, the circles to the elements in $B$  and the
squares to the e\-le\-ments in ${\rm Ap}_\mS(B)$.} \label{fig:apery2} 
\end{figure}
\end{ex}

As a direct consequence of Theorem \ref{thm:apery}, the number of elements of the Ap\'ery set ${\rm Ap}_{\mS}(B)$ coincides with the dimension
of the $\mathbb K$-vector space $\mathbb \K[\mathbf{x}]/J$. Thus, ${\rm Ap}_{\mS}(B)$ is finite if and only if $\mathbb \K[\mathbf{x}]/J$ is $0$-dimensional or, equivalently, $J \cap \mathbb K[x_i] \not= (0)$ for all $i \in \{1,\ldots,n\}$. 
The rest of this section is devoted to characterizing when this happens.

\begin{defi}
\label{def:polyhedral}
Let $\mA=\left\{ \ba_1, \ldots, \ba_n\right\}\subseteq \zz^m \oplus T$. The {\em rational polyhedral cone} $\C_\mA \subseteq \mathbb R^m$ generated by  $\mA$ is  
\[
\C_\mA = \mathrm{Cone}(\mA) \eqdef \left\{ \sum_{i=1}^n \alpha_i\, \pi(\ba_i) \mid \alpha_i \in \mathbb R_{\geq 0}\right\},
\]
where $\pi: \zz^m \oplus T \longrightarrow \zz^m$ is the canonical projection, see (\ref{canonicalProj}).

We say that $\mathcal F\subseteq \C_\mA$ is a {\em face} of $\C_\mA$ if there exists $\mathbf w\in \mathbb R^m$ such that $\mathbf w \cdot x \geq 0$ for all $x\in \C_\mA$ (where $\cdot$ represents the usual inner product) and  $\mathcal F = \{ x\in \C_\mA \mid \mathbf w \cdot x = 0 \}$.
An {\em extremal ray} of the cone $\C_\mA$ is a half-line face of $\C_\mA$.

\end{defi}

\begin{remark}
\label{remark:cone}
In the forthcoming we need the following properties of rational polyhedral cones (see, e.g.,  \cite[Proposition 1.2.12 and  Lemma 1.2.15]{cox:2011}).
\begin{enumerate}
\item \label{remark:cone-1} $\{\bz\}$ is a face of $\C_\mA$ if and only if $\pi(\mS)$ is reduced, where $\mS=\langle \mA \rangle$. 
\item  \label{remark:cone-2} Given a set $B=\left\{ \bb_1, \ldots, \bb_s\right\} \subseteq \mS \setminus \{\bz\}$ with $\mS=\langle \mA \rangle$. Then, $\C_\mA = \C_B$ if and only if for each extremal ray $r$ of $\C_\mA$, there exists $i \in \{1,\ldots,s\}$ such that $\pi(\bb_i)\in r$.
\end{enumerate}
\end{remark}

We observe that $\mS \subseteq \mathbb Z^m \oplus T$ is a reduced monoid if and only if $\pi(\mS)\subseteq \mathbb Z^m$ is reduced and $\mS \cap T =  \{\bz\}$. Thus, by the first part of Remark \ref{remark:cone}, whenever $\mS$ is a reduced monoid, then $\{\bz\}$ is a face of $\C_\mA$.

Before proceeding with the characterization of the finiteness of the Ap\'ery set $\mathrm{Ap}_{\mS}(B)$, we need a lemma in which the reduced condition of the monoid plays an important role.

\begin{lema}
\label{lemma:apfinito}
Let $\mS= \langle \ba_1, \ldots, \ba_n \rangle\subseteq \zz^m \oplus T$ be a reduced monoid and $B=\{ \bb_1, \ldots, \bb_s\} \subseteq \mS \setminus \{\bz\}$. Then, $x\in \mS$ if and only if there exist $\lambda_1, \ldots, \lambda_s \in \n$ such that 
$x-\lambda_1 \bb_1 - \cdots - \lambda_s \bb_s \in \mathrm{Ap}_\mS(B)$.
\end{lema}

\proof
Since $\mathrm{Ap}_\mS(B) \subseteq \mS$ and $B \subseteq \mS$, the claim is evident in one direction. So assume that $x\in \mS$, we will prove that 
there exist $\lambda_1, \ldots, \lambda_s \in \n$ such that 
\begin{equation}
\label{lemma:induction}
x-\sum_{i=1}^s \lambda_i \bb_i \in \mathrm{Ap}_{\mS}(B).
\end{equation}
By the first part of Remark \ref{remark:cone}, since $\mS$ is reduced, then $\{\mathbf 0\}$ is a face of $\C_{\mA}$. Therefore, there exists $\mathbf w \in \mathbb Z^n$ such that $\mathbf w \cdot \pi(x) \geq 0$ for all $x\in \mS$ and if $\mathbf w \cdot \pi(x) = 0$, then $\pi(x)=\mathbf 0$. Now we prove the lemma by induction on the value $\mathbf w \cdot \pi(x)\in \n$. If $\mathbf w\cdot \pi(x) = 0$, then $\pi(x) = \mathbf 0 \in \mathbb Z^m$, and we get that $x=\mathbf 0$ because $\mS$ is reduced. Hence, 
$x=\mathbf 0\in \mathrm{Ap}_{\mS}(B)$ and the result is true for $\lambda_1=\cdots = \lambda_s=0$. Assuming \eqref{lemma:induction} holds for any $\tilde{x}\in \mS$ such that $\mathbf w\cdot \pi(\tilde{x}) < \alpha$, for some positive integer  $\alpha$, we will prove the statement for $x\in \mS$ with $\mathbf w \cdot \pi(x) = \alpha$. We distinguish two cases: if $x\in \mathrm{Ap}_{\mS}(B)$, then it suffices to take $\lambda_1= \cdots = \lambda_s = 0$. Otherwise, by definition of the Ap\'ery set there exists $i\in \{1, \ldots, s\}$ such that $x-\bb_i \in \mS$. Let $\tilde{x} = x-\bb_i$. Then
$$\mathbf w \cdot \pi(x) = \mathbf w \cdot \pi(\bb_i) + \mathbf w \cdot \pi(\tilde{x}) \hbox{ with } \mathbf w \cdot \pi(\bb_i) > 0.$$
Thus, $\mathbf w\cdot \pi(x) > \mathbf w \cdot \pi(\tilde{x})$. We conclude, by the principle of induction, that there exist $\beta_1, \ldots, \beta_s \in \n$ such that $\tilde{x}= \sum_{j=1}^s \beta_j \bb_j \in \mathrm{Ap}_{\mS}(B)$, hence 
\[
x-\bb_i  - \sum_{j=1}^s \beta_j \bb_j \in \mathrm{Ap}_{\mS}(B).
\]
We finish the proof putting $\lambda_{j}=\beta_{j}$ for $j\in \{1, \ldots, s\}$, $j\neq i$, and $\lambda_{i}=\beta_{i}+1.$

\eproof

Let $I$ be a nonempty subset of an abelian monoid $\mS$, we say that $I$ is an {\it ideal} of  $\mS$, if for every $x \in I$ we have  $x + \mS \subseteq I$. An ideal $I \subseteq \mS$ is finitely generated if there exists a finite set $B = \{\bb_1,\ldots,\bb_s\}$ such that $I = \cup_{i = 1}^s (\bb_i + \mS)$. Clearly, in this setting we have that $$\mS \setminus I = \bigcap_{i = 1}^s\ {\rm Ap}_\mS(\{ \bb_i\}) =  {\rm Ap}_\mS(B).$$
Thus, the complement of ${\rm Ap}_\mS(B)$ in $\mS$ is just the ideal of $\mS$ spanned by $B$.

Now we can proceed with the desired characterization. Interestingly, this result also provides a criterion to determine when $I \cup \{\bz\}$ inherits  the reduced monoid structure of $\mS$, being $I$ a finitely generated ideal of~$\mS$.

\begin{thm} \label{teo:apfinito}
Let $\mS=\langle \mA \rangle = \langle \ba_1, \ldots, \ba_n \rangle\subseteq \zz^m \oplus T$ be a reduced monoid, $B=\{ \bb_1, \ldots, \bb_s\} \subseteq \mS \setminus \{\bz\}$ and  $I= \cup_{i=1}^s (\bb_i + \mS)$. The following statements are equivalent:
\begin{itemize}
\item[$(1)$] The Ap\'ery set $\mathrm{Ap}_\mS(B)$ is finite.
\item[$(2)$] $\C_\mA = \C_B$.
\item[$(3)$] $I\cup \{\mathbf 0\}$ is a (finitely generated) reduced monoid.
\end{itemize}
\end{thm}

\proof 

\noindent $(1) \Longrightarrow (3)$ Since $I \subseteq \mS$ then $I \cup \{\bz\}$ is reduced, so we just have to prove
that it is a finitely generated monoid. Assuming that  ${\rm Ap}_{\mS}(B) = \{ \mathbf h_1 = \bz, \mathbf h_2, \ldots, \mathbf h_l\}$ we will prove that
$$I \cup \{ \mathbf 0\} = \langle \{ \mathbf h_i + \mathbf b_j \mid 1 \leq i \leq l \hbox{ and } 1 \leq j \leq s\}\rangle.$$

Let $x \in I$, using Lemma \ref{lemma:apfinito}, there exist $\lambda_1, \ldots, \lambda_s \in \n$
in such a way that $x- \sum_{j=1}^s \lambda_j \bb_j \in {\rm Ap}_{\mS} (B)$. That is, there exists $i \in \{1, \ldots, s\}$ such that $\mathbf h_i = x- \sum_{j=1}^s \lambda_j \bb_j$, where not all $\lambda_j$'s are zero, since $x\notin \mathrm{Ap}_{\mS}(B)$. Thus, without loss of generality, one can assume that $\lambda_1 \neq 0$ and we can write
$$x = \mathbf h_i + \sum_{j=1}^s \lambda_j \bb_j = \mathbf h_i + \bb_1 + (\lambda_1 -1) \bb_1 + \sum_{j=2}^s \lambda_j \bb_j.$$
Hence, $x$ belongs to $\langle \{ \mathbf h_i + \mathbf b_j \mid 1 \leq i \leq l \hbox{ and } 1 \leq j \leq s\}\rangle$. The other inclusion is evident.

\noindent $(3) \Longrightarrow (2)$
In this part we are using that the unique minimal system of generators of a reduced monoid $J \subseteq \zz^m \oplus T$ consists of its irreducible elements, which is
\begin{equation}
\label{eq-pointedSemigroup}
J^{\star} \setminus  (J^{\star} + J^{\star}),
\end{equation}
where $J^{\star} = J \setminus \{\mathbf 0\}$. 

Suppose, contrary to our claim and using Remark \ref{remark:cone}(\ref{remark:cone-2}), that $\C_{\mA} \neq \C_{B}$. Then there exists an extremal ray $r$ of the cone $\C_{\mA}$ such that $\pi(\bb_i) \notin r$, for all $i \in \{ 1, \ldots, s\}$. By Definition \ref{def:polyhedral}, there exists $\mathbf w\in \mathbb R^m$ such that 
\begin{eqnarray*}
\mathbf w \cdot x \geq 0 \hbox{ for all }x\in \C_{\mA}
\hbox{, and if } x \in \C_{\mA}, \hbox{ then } \mathbf w \cdot x = 0 \Longleftrightarrow x \in r.
\end{eqnarray*} 
We define 
$\delta = \min \{ \mathbf w \cdot \pi(\bb_i) \mid 1\leq i \leq s\}$. Note that $\delta > 0$, since $\pi(\bb_i) \notin r$ for all $i \in \{1, \ldots, s\}$. We can deduce the following statements:
\begin{itemize}
\item[(a)] If $\bb\in I$ and $\mathbf w \cdot \pi(\bb) = \delta$, then we claim that $\bb \notin I+I$ and we can conclude by \eqref{eq-pointedSemigroup} that $\bb$ belongs to the minimal system of generators of $I\cup \{\bf 0\}$. Indeed, if $\bb \in I+I$ then, we can write $\bb = \bb_i + \mathbf s_1 + \bb_j + \mathbf s_2$, with $\mathbf s_1, \mathbf s_2 \in \mS$ and $i,j \in \{1, \ldots, s\}$. Hence
$$\mathbf w \cdot \pi(\bb) = \mathbf w \cdot \pi(\bb_i) + \mathbf w \cdot \pi(\mathbf s_1) + \mathbf w \cdot \pi(\bb_j) + \mathbf w \cdot \mathbf \pi(s_2) \geq 2\delta > \delta,$$
which is a contradiction.
\item[(b)] If we take $\bb_i$ such that $\mathbf w \cdot \pi(\bb_i) = \delta$ and $\ba_j \in r$, then $\mathbf w \cdot \pi(\bb_i + \lambda \ba_j) = \delta$ for all $\lambda\in \n$.
\end{itemize}
Using (a) and (b) we have actually showed that the minimal system of generators of $I\cup \{\bf 0\}$ is infinite, which contradicts our assumption.

\noindent $(2) \Longrightarrow (1)$. By Theorem \ref{thm:apery}, in order to prove that $\mathrm{Ap}_\mS(B)$ is finite it suffices to show that $\mathbb \K[\mathbf{x}] / (I_\mS + \langle \mathbf{x}^{\bbeta_1}, \ldots, \mathbf{x}^{\bbeta_s}\rangle)$ is a finite dimensional $\K$-vector space. Equivalently, we will show that there exists $g_i \in \mathbb K[x_i]$ such that $g_i(x) \in I_\mS + \langle \mathbf{x}^{\bbeta_1}, \ldots, \mathbf{x}^{\bbeta_s} \rangle$ for all $i\in \{1, \ldots, n\}$. 
In fact, we will see that there exists $\gamma_i \in \zz^+$ such that $x_i^{\gamma_i} \in I_\mS + \langle \bx^{\bbeta_1}, \ldots, \bx^{\bbeta_s} \rangle$ for all $i\in \{1, \ldots, n\}$.

Since $\C_\mA = \C_B$ and $\pi(\ba_i)\in \C_\mA$, then 
$\pi(\ba_i) = \sum_{j=1}^s \nu_j \pi(\bb_j)$ with $\nu_1, \ldots, \nu_s\in \mathbb Q_{\geq 0}.$
Thus, multiplying by an adequate positive integer $\nu$ we deduce that 
\[
\nu\, \pi(\ba_i) = \sum_{j=1}^s  \delta_j\,  \pi(\bb_j) \in \zz^m, \hbox{ where the }\delta_j \in \n \hbox{ are not all zero}.
\]
Now, multiplying by $t$, the order of $T$ we get that
\[
t \nu \ba_i = \sum_{j=1}^s  t \delta_j  \bb_j \in \zz^m \oplus T.
\]
 Hence, $x_i^{t\nu} - \prod_{j=1}^s(\mathbf{x}^{\bbeta_j})^{t\delta_j}\in I_\mS$ and we conclude that $x_i^{t \nu}\in I_\mS + \langle \mathbf{x}^{\bbeta_1}, \ldots, \mathbf{x}^{\bbeta_s} \rangle$.

\eproof

In \cite[Lemma 1.2]{pison:2003} Pis\'on gives other equivalent condition in terms of  Gr\"obner basis. However we put in value here that our proof is free of Gr\"obner bases.

\section{Elements in a reduced monoid with factorizations of the same length}
\label{sec:samelenght}

Let $\mS=\left\langle \ba_1, \ldots , \ba_n \right\rangle\subseteq \zz^m \oplus T$ be a reduced monoid given by its minimal set of generators. We consider the following subsets of $\mS$:
\[
\mL_{\mS} =  \left\{ \bb \in \mS \mid \bb \hbox{ has (at least) two different factorizations of the same length}\right\},
\]

and
\[
\mT_{\mS} = \left\{ \bb \in \mS \mid \bb \hbox{ has (at least) two different factorizations}\right\}.
\]

Observe that if $\mL_{\mS}\neq \emptyset$  (respectively  $\mT_{\mS}\neq \emptyset$) and $\bb \in \mL_{\mS}$ 
(respectively in $\mT_{\mS}$) and  $\bc \in \mS$, then $\bb + \bc \in \mL_{\mS}$ (respectively  $\mT_{\mS}$). Hence if $\mL_{\mS} \neq \emptyset$ then $\mL_{\mS}$ is an ideal of $\mS$. \\

%
%

 The next proposition shows how to obtain the set $\mT_{\mS}$ from a set
 of $\mS$-homogeneous generators of $I_\mS$. Since $I_\mS$ is a binomial ideal one may consider binomial generating sets of $I_\mS$; indeed,
 all its reduced Gr\"obner bases consist of binomials.
 
 \begin{prop}\label{prop:dosescrituras}
Let $\mS\subseteq \mathbb Z^m \oplus T$ be a reduced monoid. We get
\begin{enumerate}
\item  $\mT_{\mS} = \emptyset$ if and only if $I_{\mS} = \langle 0 \rangle$. 
\item If $I_\mS \neq \langle 0 \rangle$ and $\{g_1, \ldots , g_s\}$ is a binomial generating set of $I_\mS$ then,
$$\mT_{\mS}= \left(\deg_\mS(g_1) + \mS\right) \cup \cdots \cup \left( \deg_\mS(g_s)+\mS\right).$$
\end{enumerate}
\end{prop} 
\proof
By \eqref{eq:binomios}, we have that $b \in \mT_{\mS}$ if and only if there exists 
 a binomial $f \in I_\mS$ with ${\rm deg}_\mS(f) = b$.
 
  Since $g_i$ is a binomial in $I_\mS$, then it is $\mS$-homogeneous and $\deg_\mS(g_i) \in \mT_{\mS}$. Considering that $\mT_{\mS}$ is an ideal of $\mS$, one inclusion holds. To prove the converse, let $\bb\in \mT_{\mS}$, then there exists $f = \mathbf{x}^{\blambda}-\mathbf{x}^{\bnu}\in I_\mS$ with $\deg_\mS(\mathbf{x}^{\blambda}) = \deg_\mS(\mathbf{x}^{\bnu}) =\bb$. Now, since $I_\mS = \left\langle g_1, \ldots , g_s\right\rangle$ with $g_i = \mathbf{x}^{\balpha_i} - \mathbf{x}^{\bbeta_i}$, for some $\balpha_i, \bbeta_i \in \n^{m}$; and $f\in I_\mS$ then, one term of one of the binomials $g_i$ divides  $\mathbf{x}^{\blambda}$. That is, $\mathbf{x}^{\blambda} = \mathbf{x}^{\balpha_i}\mathbf{x}^{\bgamma}$ or equivalently, $\blambda = \balpha_i + \bgamma$ for some $\bgamma \in \n^n$ and some $i\in \{1, \ldots, s\}$. Thus,
\[
\bb=\deg_\mS(\mathbf{x}^{\blambda}) = \deg_\mS(\mathbf{x}^{\balpha_i} \mathbf{x}^{\bgamma}) = \deg_\mS(g_i) + \underbrace{\gamma_1 \ba_1 + \cdots + \gamma_n \ba_n}_{s\in \mS},
\]
where $\bgamma=(\gamma_1,\ldots, \gamma_n)$.
\eproof

\medskip

One clearly has that $\mL_{\mS} \subseteq \mT_{\mS}$.  In Lemma  \ref{lemma-Ds=D} we will
obtain $\mL_{\mS}$ by means of $\mT_{\tilde{\mS}}$ for the reduced monoid $\tilde{\mS}= \left\langle (\ba_1, 1) , (\ba_2, 1), \ldots , (\ba_n,1)\right\rangle \subseteq \zz^{m+1} \oplus T$ introduced in \eqref{eq:Stilde}.
Note that $\{(\ba_1, 1) , (\ba_2, 1), \ldots , (\ba_n,1)\}$ is the minimal set of generators of $\tilde{\mS}$. 

The idea behind considering the monoid $\tilde{\mS}$ comes from the fact that the lattice ideal $I_{\tilde{\mS}}$ is generated by
the homogeneous binomials in $I_\mS$ (see, e.g., Remark \ref{rem:binomios}). Moreover,  we will exploit the fact that factorizations of the same length of an element in $\mS$ correspond to homogeneous binomials in $I_\mS$ and, thus, to binomials in $I_{\tilde{\mS}}$. These ideas, in the particular context of numerical semigroups, have been extensively used in the study of the shifted family of a numerical semigroup (see, e.g., \cite{Vu:2014,CS:2016}).

\begin{lema}
\label{lemma-Ds=D}
Let $\mS=\left\langle \ba_1, \ldots, \ba_n\right\rangle \subseteq \zz^m \oplus T$ be a reduced monoid and $\tilde{\mS}$ the monoid defined as \eqref{eq:Stilde}. Then,
\[\mL_{\mS} = \left\{ \bx_1 \in \zz^m \mid 
(\bx_1, x_2) \in \mathcal T_{\tilde{\mS}} \hbox{ for some } x_2 \in \n \right\}.\]
\end{lema}

\proof
Let $\bx\in \mL_{\mS}$. There exist $\blambda, \bbeta \in \n^n$ such that
\[
\bx= \lambda_1 \ba_1 + \cdots + \lambda_n \ba_n = \beta_1 \ba_1 + \cdots + \beta_n \ba_n \hbox{ with } \ell(\blambda) = \ell(\bbeta) =: s\in \n.
\]

Thus, $(\bx,s) \in \zz^{m+1}$ and $(\bx,s)=\lambda_1 (\ba_1, 1) + \cdots + \lambda_n (\ba_n,1)=\beta_1 (\ba_1, 1) + \cdots + \beta_n (\ba_n,1)$,
or equivalently, $(\bx,s) \in \mT_{\tilde{\mS}}$. The other inclusion may be handled in the same way.
\eproof

The following proposition allows us to obtain $\mL_{\mS}$ from the degrees of a set of generators of the ideal $I_{\tilde{\mS}}$.

\begin{prop}
\label{prop-main-dS}
Let $\tilde{\mS}\subseteq \zz^{m+1} \oplus T$ be the monoid associated to $\mS\subseteq \zz^m \oplus T$ defined as \eqref{eq:Stilde}. We get
\begin{enumerate}
\item $\mL_{\mS} = \emptyset$ if and only if $I_{\tilde{\mS}} = \langle 0 \rangle$. 
\item If $I_{\tilde{\mS}} \neq \langle 0 \rangle$ and $\{g_1, \ldots g_s\}$ is a binomial generating set of $I_{\tilde{\mS}}$, then,
\[
\mL_{\mS} = \left( \deg_\mS(g_1) + \mS\right) \cup \cdots \cup \left( \deg_\mS(g_s) +\mS  \right).
\]
\end{enumerate}
\end{prop}

%
\proof By Lemma \ref{lemma-Ds=D} we have  $\mL_\mS = p(\mT_{\tilde{\mS}})$, where $\begin{array}{cccc} p: &\zz^{m+1} \oplus T & \rightarrow & \zz^m \oplus T \end{array}$ denotes the canonical projection. The result follows applying Proposition~\ref{prop:dosescrituras} to $\tilde{\mS}$ and observing that   $p(\deg_{\tilde{\mS}}(h)) = \deg_{\mS}(h)$ for every binomial $h \in I_{\tilde{\mS}}$.
\eproof


As a consequence of Proposition \ref{prop-main-dS}, we get the main result of this section. This result describes the set $\mS \setminus \mL_\mS$
as a particular Ap\'ery set of $\mS$.

\begin{thm}\label{thm:main}
Let $\mS\subseteq \zz^m \oplus T$ be a reduced monoid and $\{g_1,\ldots,g_s\}$ a binomial generating set of
$I_{\tilde{\mS}}$. Consider $B = \{\bb_1,\ldots,\bb_s\}$ with $ \bb_i := \deg_\mS(g_i)$ for all $i \in \{1,\ldots,s\}$. 
 Then,
\[ \mS \setminus \mL_{\mS} = {\rm Ap}_\mS(B). \]
\end{thm} 
\proof By Proposition \ref{prop-main-dS} we have that $\mL_{\mS} = \bigcup_{i = 1}^s (\bb_i + \mS)$. Therefore
\[ \mS \setminus \mL_{\mS} = \bigcap_{i = 1}^s\ {\rm Ap}_\mS(\{ \bb_i\}) =  {\rm Ap}_\mS(B).\]
\eproof


Theorems \ref{thm:main} and \ref{thm:apery} provide a method to compute $\mS \setminus \mL_{\mS}$.  More precisely, 
\begin{enumerate}
\item [(I)] Consider a binomial generating set $\{g_1,\ldots,g_s\}$ of $I_{\tilde{\mS}}$ and denote $g_i = \bx^{\balpha_i}  - \bx^{\bbeta_i}$ for all $i \in \{1,\ldots,s\}$ (see Remark \ref{rem:binomios}).
\item [(II)] Then, one can apply Theorem \ref{thm:apery} to compute ${\rm Ap}_\mS(B)$ being $B = \{\bb_1,\ldots,\bb_s\}$ with $\bb_i = \deg_{\mS}(g_i)$. 
\end{enumerate}
In order to use Theorem \ref{thm:apery}, as it is stated, one needs a factorization of $\bb_1,\ldots,\bb_s$. Nevertheless, this does not involve any extra computations. Indeed, $I_{\tilde{\mS}}$ is an $\mS$-homogeneous ideal and, hence,  $\balpha_i$ and $\bbeta_i$ are two factorizations of $\bb_i$ for all $i \in \{1,\ldots,s\}$. \\
Let us illustrate this method in the next example.

\begin{ex}\label{ex:aperycontinuacion} Consider, as in Example \ref{ex:apery}, the affine monoid \[ \mS = \langle \ba_1,\ldots,\ba_5\rangle \subseteq \zz^2,\] with $\ba_1 = (0,2), \ba_2 = (1,2), \ba_3 = (1,1), \ba_4 = (3,2), \ba_5 = (4,2)$ and let us compute $\mL_\mS$ and $\mS \setminus \mL_\mS$. For this purpose, we first consider $I_{\tilde{\mS}}$ with $\tilde{\mS} = \langle (0,2,1), (1,2,1),  (1,1,1), (3,2,1),  (4,2,1)\rangle$. It turns out that $I_{\tilde{\mS}}$ is minimally generated by $\{g_1,\,g_2,\,g_3,\, g_4\}$, where:
\[ \begin{array}{llll}  
g_1 = x_2^3-x_1^2 x_4, & g_2 = x_2 x_4-x_1 x_5, & g_3 = x_4^3-x_2 x_5^2, & 
g_4 = x_1 x_4^2-x_2^2 x_5.
\end{array} \]
Let $B = \{\bb_1, \bb_2, \bb_3, \bb_4\}$ where $\bb_i := \deg_\mS(g_i)$. One gets that $\bb_1 = 3 \ba_2 =(3,6),\, \bb_2 = \ba_2 + \ba_4 = (4,4), \, \bb_3 = 3 \ba_4 = (9,6)\text{ and } \bb_4 = \ba_1 + 2 \ba_4 = (6,6)$. By Proposition \ref{prop-main-dS} we have:
\[ \mL_\mS = \cup_{i = 1}^4 (\bb_i + \mS) = ((3,6) + \mS) \cup ((4,4) + \mS) \cup ((9,6) + \mS) \cup ((6,6) + \mS). \]
Moreover, since $\bb_4 = (4,4) + (2,2) \in \bb_2 + \mS$, we put 
\[ \mL_\mS = \cup_{i = 1}^3 (\bb_i + \mS) = ((3,6) + \mS) \cup ((4,4) + \mS) \cup ((9,6) + \mS). \]
Thus, setting $B' = \{\bb_1,\bb_2,\bb_3\}$ we have $\mS \setminus \mL_\mS = {\rm Ap}_\mS(B')$ and this set equals the one we computed in 
Example \ref{ex:apery}. So the squared grid points in Figure \ref{fig:apery} correspond to the elements of $\mS \setminus \mL_\mS$.
\end{ex}

\begin{ex}\label{ex:aperycontinuacion} Consider, as in Example \ref{ex:apery2}, the reduced monoid \[ \mS = \langle \ba_1,\ba_2,\ba_3\rangle \subseteq \zz \oplus \zz_2,\] with $\ba_1 = (2, \overline{0}), \ba_2 = (3, \overline{1}), \ba_3 = (4, \overline{1})$ and let us compute $\mL_\mS$ and $\mS \setminus \mL_\mS$. For this purpose, we first consider $I_{\tilde{\mS}}$ with $\tilde{\mS} = \langle (2,1,\overline{0}), (3,1, \overline{1}),  (4,1, \overline{1})\rangle \subseteq \mathbb Z^2 \oplus \mathbb Z_2$. It turns out that $I_{\tilde{\mS}} = \langle g \rangle$ with $g=x_2^4-x_1^2x_3^2$.
Let $B = \{\bb_1\}$, where $\bb_1 = \deg_{\mS}(g) = 4 \cdot \ba_2 = (12, \overline{0})$, by Proposition \ref{prop-main-dS} we have:
\[ \mL_\mS = (\bb_1 + \mS) = ((12, \overline{0}) + \mS). \]
Thus, we have $\mS \setminus \mL_\mS = {\rm Ap}_\mS(B)$ and this set equals the one we computed in 
Example \ref{ex:apery2}. 
\end{ex}

As a direct consequence of Theorems \ref{thm:apery} and \ref{thm:main}, we have:

\begin{cor}\label{cor:dimen}
Let $\mS\subseteq \zz^m \oplus T$, be a reduced monoid. Then, 
\[\sharp(\mS\setminus\mL_{\mS}) = \dim \left( \mathbb \K[\mathbf{x}] / (I_\mS + \mathrm{in}_{\succ} (I_{\tilde{\mS}}))\right),\]
where $\mathrm{in}_{\succ} (I_{\tilde{\mS}})$ represents the initial ideal of $I_{\tilde{\mS}}$ with respect to any monomial order.
\end{cor}

Now, putting this result together with Theorem  \ref{teo:apfinito} we get the following corollary, characterizing when there is only a finite number of elements of $\mS$ not belonging to $\mL_{\mS}$. It is also worth mentioning that this happens if and only if $\mL_\mS \cup \{\bz\}$ inherits the finitely generated reduced monoid structure of $\mS$. 

\begin{cor}
\label{cor:maincoro}
Let $\mS = \langle \mA \rangle \subseteq \zz^m \oplus T$ be a finitely generated reduced monoid. Then, the following statements are equivalent:
\begin{itemize}
\item[{\rm (1)}] $\mS\setminus\mL_{\mS}$ is a finite set.
\item[{\rm (2)}] For every extremal ray $r$ of $\C_\mA$ there are either:
\begin{itemize} \item[{\rm (2.a)}] two elements $\ba_1,\ba_2 \in \mA$ such that $\pi(\ba_1) = \pi(\ba_2) \in r$, or 
\item[{\rm (2.b)}] three elements $\ba_1,\ba_2,\ba_3 \in \mA$ such that $\pi(\ba_1),\pi(\ba_2),\pi(\ba_3) \in r$.
\end{itemize}
\item[{\rm (3)}] $\mL_{\mS} \cup \{\mathbf{0}\}$ is (a finitely generated) reduced monoid.
\end{itemize}
\end{cor} 
\proof  Being (1) and (3) equivalent by Theorem \ref{teo:apfinito}, we are going to prove the equivalence between (1) and (2). Let $I_{\tilde{\mS}} = (g_1,\ldots,g_s),$ where $g_i$ is a binomial and $B =\{\bb_1,\ldots,\bb_s\}$, with $\bb_i := \deg_\mS(g_i)$.
 By  Theorem \ref{thm:main}  we have  $\mS \setminus \mL_{\mS} = {\rm Ap}_\mS(B)$. Thus, by Proposition \ref{teo:apfinito} and Remark \ref{remark:cone},  $\mS \setminus \mL_{\mS}$ is finite if and only if there is at least one element of $\pi(B)$ in each extremal ray of $\C_\mA$. So it just remains to prove that this happens if and only if either (2.a) or (2.b) holds. Consider an extremal ray $r$. We take $R = \langle (\ba_1,1), (\ba_2,1) \rangle$ if (2.a) holds, or $R = \langle (\ba_1,1), (\ba_2,1), (\ba_3,1)\rangle$ if (2.b) holds. In both cases we have that $I_R$ is a height one lattice ideal (see Remark \ref{rem:binomios}). Then, there is a binomial $f \in I_R \subseteq I_{\tilde{\mS}}$.  As a consequence, one of the monomials appearing in $g_1,\ldots,g_s$ has to divide one of the monomials appearing in $f$. Hence, the $\mS$-degree of the corresponding $g_i$ belongs to $R$ and, $\pi(\bb_i) = \pi({\rm deg}_\mS(g_i)) \in r$.  Conversely, if $\pi(\bb_i)$ is in $r$, then we have $g_i = \mathbf{x}^{\balpha_i} - \mathbf{x}^{\bbeta_i} \in I_{\tilde{\mS}}$ and we may assume that $\mathbf{x}^{\balpha_i} $ and $\mathbf{x}^{\bbeta_i}$ are relatively prime. Since $g_i$ is homogeneous, then:
\begin{itemize}
\item[{\rm (a)}] either $g_i = x_j^d - x_k^d$ with $d \ba_j = d \ba_k$, 
\item[{\rm (b)}] or there are at least three variables involved in $g_i$. 
\end{itemize}
If (a) holds, then $d \pi(\ba_j) = d \pi(\ba_k) \in r$, and $\pi(\ba_j) = \pi(\ba_k) \in r$.  If (b) holds, given that $\pi(\bb_i) = \sum_{i = 1}^n \alpha_{ij} \pi(\ba_j) = \sum_{i = 1}^n \beta_{ij} \pi(\ba_j)$ and $r$ is an extremal ray, we have  $\pi(\ba_j) \in r$ whenever $\alpha_{ij} \neq 0$ or $\beta_{ij} \neq 0$. Hence there are at least three $\pi(\ba_i)$ in $r$, finishing the proof.
\eproof 

\medskip

Observe that condition (2.a) cannot occur when $\mS \subseteq \zz^m$ is an affine monoid.

In the remainder of the section we will apply our study to the setting of numerical semigroups. More precisely, we will deduce the results of \cite{chapman:2011}, using Proposition \ref{prop-main-dS}, in the setting of numerical semigroups.

Let $\mS = \langle a_1,\ldots, a_n\rangle \subseteq \n$ be a numerical semigroup given by its minimal generating set.
Denote by $F(\mS)$ the Frobenius number of $\mS$, which is the largest integer not in $\mS$, i.e., $F(\mS) = \max (\zz\setminus \mS)$.

We reprove \cite[Proposition 2]{chapman:2011} in the next corollary.
\begin{cor}\label{cor:chapman}
Let $\mS = \langle a_1,\ldots, a_n\rangle \subseteq \n$ be a numerical semigroup with Frobenius number  $F(\mS)$.
\begin{enumerate}
\item Let $w\in \mL_{\mS}$. For any integer  $z$ verifying $z > w+F(\mS)$, we have  $z\in \mL_{\mS}$.
\item If $\mL_{\mS}\neq \emptyset$ then $\mL_{\mS} \cup \{0\}$ is a numerical semigroup.
\end{enumerate}
\end{cor}
\proof
Remember that  $\mL_{\mS}$ is a semigroup. The first assertion follows from the definition of $F(\mS)$.
On the other hand, since $\mL_{\mS}$ is a semigroup then  $\mL_{\mS} \cup \{0\}$ is a submonoid of $\n$. By the first assertion of this corollary, $\mL_{\mS} \cup \{0\}$ has finite complement in $\n$.
\eproof

%
%
%

Now, we reprove  \cite[Theorems 2 and 3]{chapman:2011}:
\begin{cor} \label{cor:numericos}
Let $\mS=\left\langle a_1, \ldots , a_n \right\rangle\subseteq \n$ be a numerical semigroup given by its minimal set of generators.
\begin{enumerate}
\item $\mL_{\mS}= \emptyset$ if and only if $n \leq 2$. 
\item If $n=3$, then $\mL_{\mS}=  (a_2 (a_3-a_1) / \gcd(a_2-a_1,a_3-a_1) ) + \mS$.
\end{enumerate}
\end{cor}

\proof
The height of the ideal $I_{\tilde{\mS}}$ equals  $\max\{0,n-2\}$ (see Remark \ref{rem:binomios}). Thus, $I_{\tilde{\mS}} = \langle 0 \rangle$ if and only if $n \leq 2$. If $n = 3$, then $I_{\tilde{\mS}}$ is the principal ideal 
\[
    I_{\tilde{\mS}} = \left\langle x_2^{(a_3-a_1)/d}- x_1^{(a_3 - a_2)/d} x_3^{(a_2 - a_1)/d}\right\rangle, 
\]
with $d := \gcd(a_2-a_1,a_3-a_1)$. Thus, by Proposition \ref{prop-main-dS}, we conclude  \[ \mL_{\mS}= \deg_\mS\left(x_2^{(a_3-a_1)/d}\right)  +\mS =  \left(a_2(a_3-a_1)/d  \right) + \mS. \]
\eproof

\section{The equal catenary degree}
\label{sec:catenaryDegree}

Let $\mS = \left\langle \ba_1, \ldots, \ba_n \right\rangle \subseteq \mathbb Z^m \oplus T$ be a reduced monoid given by its minimal set of generators. Let $\blambda = (\lambda_1, \ldots, \lambda_n) \in \mathbb N^n$ and $\bnu = (\nu_1, \ldots, \nu_n)\in \mathbb N^n$ be two factorizations of the same length of an element $\bb \in \mS$. We define the {\it distance} between $\blambda$ and $\bnu$ as:
\[d(\blambda, \bnu) = \sum_{i=1}^n \left(\lambda_i -  \min \{\lambda_i ,\nu_i \}\right) = \sum_{i=1}^n \left(\nu_i -  \min \{\lambda_i ,\nu_i \}\right) .\]

Let $N \in \mathbb N$, a finite sequence $(\blambda = \bgamma_0 , \bgamma_1, \ldots, \bgamma_k = \bnu)$ of factorizations of $\bb \in \mS$ of the same length is called an $N${\it -chain}  from $\blambda$ to $\bnu$ if $d(\bgamma_{i-1}, \bgamma_i) \leq N$ for all $i \in \{1,\ldots, k\}$. In what follows, when we say an $N$-chain we mean an $N$-chain of factorizations of the same length.

Let $\rm c_{\rm eq}(\bb)$ denote the smallest $N \in \mathbb N \cup \{\infty\}$ with the following property: for any $\blambda$, $\bnu$ factorizations of $\bb$ of the same length, there exists an $N$-chain from $\blambda$ to $\bnu$. That is,
\[{\rm c_{\rm eq}}(\bb) = \min \left\{ N \in \mathbb N \cup \{\infty\} \mid 
\begin{array}{c}
\hbox{there exists an } N\hbox{-chain for any two  }\\
  \hbox{ factorizations of the same length of } \bb
\end{array}
\right\}.\]
The value
$\rm c_{\rm eq}(\mS) = \max \{ \rm c_{\rm eq}(\bb) \mid \bb \in \mS\}$
is called the \emph{equal catenary degree of}~$\mS$.

 Equal catenary degrees have been studied since $2006$, see for example \cite{foroutan:2006,geroldinger:2006,hassler:2009,blanco:2011,philipp:2015,geroldinger:2019} and the references therein.\\

From the definition it follows that an element $\bb \in \mS$ has equal catenary degree $\rm c_{\rm eq}(\bb) >0$ if and only if it has, at least,  two different factorizations of the same length. As a consequence, $\rm c_{\rm eq}(\mS) > 0$ if and only if $\mL_\mS \not= \emptyset$ which, by Proposition \ref{prop-main-dS}, is equivalent to $I_{\tilde{\mS}} \neq (0)$. In this section we dig deeper into the connections between $\rm c_{\rm eq}(\mS)$ and the ideal $I_{\tilde{\mS}}$. The main result in this section is Theorem \ref{th:catenary}, where we prove that $\rm c_{\rm eq}(\mS)$ equals the maximum degree of the elements of a minimal set of homogeneous generators of $I_{\tilde{\mS}}$. To prove this result we use the following remark.
  
\begin{remark}
\label{lemma:catenary}
Let $\mS = \left\langle \ba_1, \ldots, \ba_n \right\rangle \subseteq \mathbb Z^m \oplus T$ be a reduced monoid and let $\blambda, \bnu \in \mathbb N^n$ be two factorizations of $\bb\in \mS$ of the same length.
Then, $\mathbf x^{\blambda} - \mathbf x^{\bnu} = \gcd (\mathbf x^{\blambda}, \mathbf x^{\bnu}) \cdot f$, where $f \in I_{\tilde\mS}$ is a binomial of degree $\deg(f) = d(\blambda, \bnu)$.
\end{remark}

Let us proceed with Theorem \ref{th:catenary}. The inequality in the second part of  this theorem already appears in \cite[Proposition 4.4.3]{blanco:2011}. Anyway we add its proof in order to improve readability of the article.

\begin{thm}
\label{th:catenary}
Let $\mS= \left\langle \ba_1, \ldots, \ba_n\right\rangle \subseteq \mathbb Z^m \oplus T$ be a reduced monoid and set $\tilde{\mS} = \left\langle (\ba_1, 1) , \ldots, (\ba_n, 1)\right\rangle \subseteq \mathbb Z^{m+1}\oplus T$. We have that 
\begin{enumerate}
\item ${\rm c_{\rm eq}}(\mS) = 0$ if and only if $I_{\tilde{\mS}} = \langle 0 \rangle$. 
\item If $I_{\tilde{\mS}} \neq \langle 0 \rangle$ and $\{f_1, \ldots, f_s\}$ is a binomial generating set of $I_{\tilde{\mS}}$, then
\[{\rm c_{\rm eq}}(\mS) \leq \max_{1\leq i \leq s} \{ \deg(f_i)\}.\]
\item Moreover, if $\{f_1, \ldots, f_s\}$ is a binomial minimal generating set of $I_{\hat{\mS}}$, then
\[{\rm c_{\rm eq}}(\mS) = \max_{1\leq i \leq s}\{ \deg(f_i)\}.\]
\end{enumerate}
\end{thm}

\proof
The first statement follows from the definition of equal catenary degree and Proposition \ref{prop-main-dS}. So assume that $I_{\tilde{\mS}} \neq \langle 0 \rangle$. 
Let $\{f_1, \ldots, f_s\}$ be a binomial generating set of $I_{\tilde{\mS}}$. 
Put $M:= \max_{1\leq i \leq s} \{ \deg(f_i)\}$. Let us prove that ${\rm c_{\rm eq}}(\mS) \leq M$. Consider  two factorizations  of the same length, $\blambda, \bdelta \in \mathbb N^n$, of an element $\bb \in \mS$. Let us find an $M$-chain between them. 
Since $g := \bx^{\blambda} - \bx^{\bdelta} \in I_{\tilde{\mS}}$ and $\{f_1,\ldots,f_s\}$ is a binomial generating set, then $g$ can be written as (see, e.g., \cite[Proposition 3.11]{CTV2017}) \[g = \sum_{j = 1}^k \bx^{\bnu_j} h_{j},\] with $ h_j  \in \{\pm f_1,\ldots,\pm f_s\}$, and if we put $h_j = \bx^{\balpha_j} - \bx^{\bbeta_j}$ with $\balpha_j,\bbeta_j \in \mathbb{N}^n$, then $\blambda = \bnu_1 + \balpha_1,\, \bnu_i + \bbeta_i = \bnu_{i+1} \balpha_{i+1}$ for all $i \in \{1,\ldots,k-1\},$ and $\bdelta = \bnu_{k} + \bbeta_{k}$.
As a consequence \[ (\blambda = \bnu_1 + \balpha_1, \bnu_1 + \bbeta_1 = \bnu_2 + \balpha_2, \ldots, \bnu_{k-1} + \bbeta_{k-1} = \bnu_k + \balpha_k,
\bnu_{k} + \bbeta_{k} = \bdelta) \] is an $M$-chain, because $h_j$ is a homogeneous element of $I_{\tilde{\mS}}$ and $d(\bnu_j + \balpha_j, \bnu_j + \bbeta_j) = {\rm deg}(h_j) \leq M$ for all $j \in \{1,\ldots,k\}$.

Now take $\{f_1,\ldots,f_s\}$ a minimal set of generators of $I_{\hat{\mS}}$. Fix  $i \in \{1,\ldots,s\}$ and write $f_i = \mathbf x^{\balpha} - \mathbf x^{\bbeta}$. Set $M = \deg(f_i)$, then $(\balpha, \bbeta)$ is an $M$-chain from $\balpha$ to $\bbeta$. 
We claim that there is no  $N$-chain from $\balpha$ to $\bbeta$ for all $N < M$. 
 Suppose, contrary to our claim, that there exists an $N$-chain of factorizations of $\bb := \deg_{\mS}(\mathbf x^{\balpha}) = \deg_{\mS}(\mathbf x^{\bbeta})$ from $\balpha$ to $\bbeta$ with $N< M$. That is, there exists a finite sequence $(\balpha= \bgamma_0, \bgamma_1, \ldots, \bgamma_k = \bbeta)$ of factorizations of the same length of $\bb$ with $d(\bgamma_{j-1}, \bgamma_{j}) \leq N$ for all $j \in \{1, \ldots, k\}$. Thus, by Remark \ref{lemma:catenary}, there exist $g_1, \ldots, g_k \in I_{\hat{\mS}}$ such that
$$f_i = \mathbf x^{\balpha} - \mathbf x^{\bbeta} = \sum_{j = 1}^k \left(x^{\bgamma_{j-1}} - x^{\bgamma_{j}}\right) = \sum_{j=1}^k \mathbf x^{\bdelta_j} g_j,$$
with $\deg(g_j) = d(\bgamma_{j-1}, \bgamma_{j}) \leq N < M$, which contradicts the minimality of $\{ f_1, \ldots, f_s\}$. This implies that  ${\rm c_{\rm eq}}(\mS) \geq \max_{1\leq i \leq s}\{ \deg(f_i)\}$ and the result follows. 
\eproof

Thus, whenever one knows an explicit set of generators of $I_{\tilde{\mS}}$, one can compute the value ${\rm c_{\rm eq}}(\mS)$. This is the case of three generated numerical semigroups, allowing us to re-prove \cite[Lemma 6]{gonzalez2019monotone}.

\begin{cor}
Let $\mS=\left\langle a_1, a_2,a_3\right\rangle \subseteq \mathbb N$ be a numerical semigroup given by its minimal set of generators. Then,
\[{\rm c_{\rm eq}}(\mS) = \frac{a_3-a_1}{\gcd (a_2-a_1, a_3-a_1)}.\]
\end{cor}

\proof
In the proof of Corollary \ref{cor:numericos} we observed that $I_{\tilde{\mS}} = \langle g \rangle$, being
$g = x_2^{(a_3-a_1)/d}- x_1^{(a_3 - a_2)/d} x_3^{(a_2 - a_1)/d}$ with $d := \gcd(a_2-a_1,a_3-a_1)$. Hence, applying Theorem \ref{th:catenary} we get ${\rm c_{\rm eq}}(\mS) = {\rm deg}(g) = (a_3-a_1)/d$.
\eproof

\medskip

Given a homogeneous ideal $J \subseteq \K[\bx]$, the Castelnuovo-Mumford regularity of $J$, denoted ${\rm reg}(J)$, is the maximum among all the values $b_j - j$, where $b_j$ is the degree of a $j$-th syzygy in a minimal graded free resolution of $J$ (see, e.g., \cite{BayerMumford,EisGoto} for other equivalent definitions). In particular, 
${\rm reg}(J)$ provides an upper bound for the degrees of the $0$-syzygies, which correspond to the degrees in a minimal generating set of $J$. As a direct consequence of Theorem \ref{th:catenary} we have that ${\rm c_{\rm eq}}(\mS) \leq {\rm reg}(I_{\tilde{S}})$ and, thus, upper bounds for ${\rm c_{\rm eq}}(\mS)$ can be derived from upper bounds on the regularity of ${\rm reg}(I_{\tilde{S}})$. We finish the section  applying this idea in the context of numerical semigroups. In order to  provide an upper bound for the equal catenary degree of any numerical semigroup
we use the upper bound for the Castelnuovo-Mumford regularity of projective monomial curves obtained by L'vovsky:

\begin{prop}{\rm \cite[Proposition 5.5]{Lvovsky}} Let $0 = b_1 < b_2 < \cdots < b_n$ a sequence of relatively prime integers and consider $\mT = \langle (b_1,1), \ldots, (b_n,1)\rangle$, then   \[ {\rm reg}(I_\mT) \leq {\rm max}_{1 \leq i < j < n}\{b_{i+1}-b_i+b_{j+1}-b_j\}.\] 
\end{prop}

Finally we need the next remark, which will also be  useful in the remaining sections.

\begin{remark}
\label{equivalent-ideal}
 Let $\mS$ be a numerical semigroup generated by $\mathcal A=\{a_1, \ldots, a_n\}\subseteq \n$ with $a_1<\cdots <a_n$. Consider the affine monoid 
\[\tilde{\mS}  = \langle  (a_1, 1) , \ldots, (a_n,1) \rangle\subseteq \n^2,\]
associated to $\mS$.
The following operations allow us to define, from $\mathcal A$, new monoids $\mathcal T\subseteq \n^2$ determining the same (toric) ideal  $I_{\tilde{\mS}}\subseteq \mathbb K[x_1, \ldots, x_n]$. 

\begin{itemize}
\item[(1)] Subtracting to each element of $\mathcal A$ the same scalar $\lambda \leq a_1$, $\lambda \in \mathbb N$. 
Consi\-dering  $\mathcal T = \langle (a_1 - \lambda,1), \ldots, (a_n - \lambda,1)\rangle\subseteq \n^2$, then 
$I_{\mathcal T} = I_{\tilde{\mS} }$.

\item[(2)] Subtracting each element of $\mathcal A$ to the same scalar $\lambda \geq a_n$, $\lambda \in  \n$.  
Consi\-dering $\mathcal T = \langle (\lambda - a_1,1), \ldots, (\lambda - a_n ,1)\rangle\subseteq \n^2$, then
$I_{\mathcal T} = I_{\tilde{\mS}}$.

\item[(3)] Multiplying and dividing all the elements of $\mathcal A$ by the same scalar.
Consi\-dering $\lambda\in \n$ a divisor of $\gcd(a_1, \ldots, a_n)$ and $\mathcal T = \langle (\frac{a_1}{\lambda},1), \ldots, (\frac{a_n}{\lambda}, 1)\rangle \subseteq \n^2$, then, $I_{\mathcal T} = I_{\tilde{\mS} }$. A similar property can be deduced if we multiply each element of $\mathcal A$ by a constant $\lambda \in \mathbb Z^+$.
\end{itemize}
\end{remark}

\begin{thm}\label{th:cotacurvas}Let $\mS \subseteq \mathbb N$ be a numerical semigroup with minimal set of generators $a_1 < \cdots < a_n$ and $n \geq 3$. Then,
\[{\rm c_{\rm eq}}(\mS) \leq \frac{{\rm max}_{1 \leq i < j < n}\{a_{i+1}-a_i + a_{j+1}-a_j\} }{\gcd (a_2-a_1, a_3-a_1,\ldots,a_n-a_1)}.\]
\end{thm}

\proof Let $\tilde \mS = \langle (a_1,1),\ldots,(a_n,1) \rangle \subseteq \mathbb N^2$.
By Theorem \ref{th:catenary} we get ${\rm c_{\rm eq}}(\mS) \leq {\rm reg}(I_{\tilde{\mS}})$.  After  Remark \ref{equivalent-ideal}.(1) for $\lambda = 1$, and then Remark \ref{equivalent-ideal}.(3) with $d = \gcd(a_2-a_1,a_3-a_1,\ldots,a_n-a_1)$;  we have  $I_{\tilde{\mS}} = I_{\mathcal T},$ where
\[\mathcal T = \langle (b_1,1), (b_2,1), \ldots, (b_n, 1) \rangle,\]
 being $b_i = \frac{a_i-a_1}{d}$ for all $i \in \{1,\ldots,n\}$.  Applying L'vovsky's bound to $I_T$ we get \[ \begin{array}{llll} {\rm reg}(I_{\tilde{\mS}}) =  {\rm reg}(I_T) & \leq & {\rm max}_{1 \leq i < j < n}\{b_{i+1}-b_i+b_{j+1}-b_j\} = \\ & = & {\rm max}_{1 \leq i < j < n}\{a_{i+1}-a_i + a_{j+1}-a_j\} / d. \end{array}\]
\eproof

\section{Computing $\mL_{\mS}$ when $\mS$ is generated by an almost arithmetic sequence}
\label{sec:arithmetic}
In this section we will focus our attention on computing $\mL_\mS$ in the particular case of numerical semigroups generated by an almost arithmetic sequence. As a warm-up we begin with the case of arithmetic sequences. 

Let $\mS$ be a numerical semigroup generated by an arithmetic sequence of relative primes, i.e., $\mS = \langle m_1, \ldots, m_n\rangle \subseteq \n$ where $m_1 < \cdots < m_n$ is an arithmetic sequence and  $\gcd(m_1, \ldots, m_n)=1$. In other words, 
\begin{equation}
\label{eq:arithmetic-sequence}
m_i = m_1 + (i-1)e \hbox{ for some $e$ with }\gcd (m_1, e) = 1\hbox{ for all }i\in \{2, \ldots, n\}.
\end{equation}
An \emph{almost arithmetic sequence} is a sequence in which all but one of the elements form an arithmetic sequence.

\begin{prop}
\label{prop:arithmetic-sequence}
Let $\mS = \langle m_1, \ldots, m_n\rangle \subseteq \n$ be a numerical semigroup generated by an arithmetic sequence of relative primes as in equation \eqref{eq:arithmetic-sequence}. Then
$$\mL_{\mS} =  \{2m_1 + \lambda e  \mid 2 \leq \lambda \leq 2n-4\} + \mS,$$
where $e := m_2 - m_1$ is the difference of the arithmetic sequence.
\end{prop}

\proof
We define $m_i' = m_i -m_1 = (i-1) e$ and $m_i'' = \frac{m_i'}{e}$. Then, by Remarks \ref{equivalent-ideal}.(1) and \ref{equivalent-ideal}.(3), we have  $I_{\tilde{\mS}} = I_{\mT_1} = I_{\mT_2}$ with 
\begin{eqnarray*}
\mT_1  & = & \langle (0,1), (m_2', 1) , \ldots, (m_{n}',1)\rangle = \langle (0,1), (e,1), \ldots, ((n-1)e, 1) \rangle \subseteq \n^2,\\
\mT_2 & = & \langle (0,1), (m_2'', 1) , \ldots, (m_{n}'',1)\rangle = \langle (0,1) , (1,1), \ldots, (n-1,1)\rangle\subseteq \n^2.
\end{eqnarray*}
Moreover, $I_{\mT_2}$ is the defining ideal of the rational normal curve in $\mathbb P_{\mathbb K}^{\, n-1}$ of degree $n-1$. Indeed, 
$I_{\mT_2} = \langle x_ix_j - x_{i-1}x_{j+1} \mid 2 \leq i \leq j \leq n-1\rangle$.
Thus, after Proposition \ref{prop-main-dS} we obtain $\mL_{\mS}$ from the set of generators of the ideal $I_{\tilde{\mS}}$. That is,
\[ \mL_{\mS} =  \{m_i + m_j \mid 2\leq i \leq j \leq n-1\} + \mS = \{2m_1 + \lambda e \, \vert \, 2 \leq \lambda \leq 2n-4\} + \mS. \]
\eproof

In the previous result we obtained an explicit minimal set of generators of~$I_{\tilde{\mS}}$. As a consequence, we get an alternative proof of \cite[Theorem 3]{gonzalez2019monotone}:
\begin{cor}Let $\mS = \langle m_1, \ldots, m_n\rangle \subseteq \n$ be a numerical semigroup gene\-ra\-ted by an arithmetic sequence, then ${\rm c_{\rm eq}}(\mS) = 2$.
\end{cor}
%
%
%

In the rest of the section we will focus on the case of numerical semigroups generated by an {\it almost arithmetic sequence}, i.e. $\mS = \langle m_1, \ldots, m_n, b \rangle \subseteq \n$ and there exists $e\in \n$ such that
\begin{equation}
\label{eq-almostarithmetic}
m_i = m_1 + (i-1)e \hbox{ with } \gcd (m_1, e ) = 1 \; \hbox{\rm  for all }i\in \{2, \ldots, n\}.
\end{equation}
In Theorem \ref{prop:almostarithmetic-A} we will provide a description of $\mL_\mS$ in this setting. In its proof we will use the following two remarks.

\begin{remark}
\label{rem:resultado-Eva}
Let $\mathcal A = \{ m_1, \ldots, m_n\} \subseteq \n$ be an arithmetic sequence of relative primes
and consider the following affine  monoid
$$\mathcal T = \langle (0,1), (m_1,1), \ldots, (m_n,1)\rangle\subseteq \n^2.$$
By  \cite[Theorem 2.2]{bermejo:2017}, the ideal $I_\mathcal T$ is minimally generated by 

\begin{equation}
\label{Resultado-Eva}
\{x_i x_j + x_{i-1}x_{j+1} \mid 2\leq i \leq j \leq n-1\} 
 \bigcup  \{ x_1^{\alpha}x_i - x_{n-k+i} x_n^{\alpha-e}x_{n+1}^e \mid 1\leq i \leq k \},
\end{equation}
where the pair $(\alpha, k) \in \n^2$ is defined as follows: 
\begin{itemize}
\item $k$ is the only integer such that $k \equiv 1-m_n \mod (n-1)$ and $1 \leq k \leq n-1$, and
\item $\alpha = \left\lfloor \frac{m_n-1}{n-1} \right\rfloor \in \n$, where $\lfloor \cdot \rfloor$ denotes the floor function.
\end{itemize}

\end{remark}

\begin{remark} \label{rm:quitaruno}
Consider the monoids
$$\begin{array}{c}
\mathcal T_1 = \langle (0,1), (a_2, 1), \ldots, (a_n,1), (b,1)\rangle\subseteq \mathbb Z^2\\
\mathcal T_2 = \langle (0,1), (a_2,1), \ldots, (a_n,1)  \rangle\subseteq \mathbb Z^2,\\
\end{array}$$
where $a_2, \ldots, a_n, b\in \mathbb Z^+$ are relatively prime. Set $B=\gcd(a_2, \ldots, a_n)$,~if
$B\cdot b = \sum_{i=2}^n \alpha_ia_i \hbox{ for some } \alpha_i \in \n \hbox{ such that } \sum_{i=1}^n \alpha_i \leq B;$
then as a direct consequence of \cite[Lemma 2.1 and Proposition 2.2]{BG:2014}, we have

\[ I_{\mathcal T_1} = I_{\mathcal T_2} \cdot \mathbb K[x_1,\ldots,x_{n+1}] + \langle x_{n+1}^B - x_1^{B-\sum_{i=1}^n \alpha_i} \prod_{i=2}^n x_{i}^{\alpha_i}  \rangle. \]
\end{remark}

Before proceeding with the proof of the main result of this section, we setup some notation.  Let $\mA=\{m_1, \ldots, m_n, b \}$ be an almost arithmetic sequence as in \eqref{eq-almostarithmetic}. Put $M:=\max \mA$, $m:=\min \mA$, $d=\gcd(b-m_1,e)$, $\beta = \left\lfloor \frac{M-m-d}{d(n-1)} \right\rfloor$ and 
$H:= \{2m_{1}+\lambda e\;\mid \; 2\leq \lambda \leq 2n-2\}$.

\begin{thm}
\label{prop:almostarithmetic-A}
Let $\mathcal A = \{ m_1, \ldots, m_n, b\} \subseteq \n$ be an almost arithmetic sequence and consider the numerical semigroup $\mS $ generated by $\mathcal A$. 

\begin{enumerate}
\item [\rm{(I)}] Suppose that $b\in \{m,M\}$.
 \end{enumerate}
 
\begin{enumerate}
\item If $d(n-1)$ divides $M-m$, then
\[
\mathcal L_{\mS} = H \cup   \left(  (\beta+1) m_1 + \mS\right),\;\; \hbox{\rm when $b=m$}
\]
or 
\[
\mathcal L_{\mS} = H \cup   \left(  (\beta+1) m_n + \mS\right),\;\; \hbox{\rm when $b=M$}. 
\]

\item If $d(n-1)$ does not divide $M-m$, then
\[
\mathcal L_{\mS} = H \cup   \left(  \{(\beta+1) m_1, (\beta + 1) m_1 + e\} + \mS \right), \;\; \hbox{\rm when $b=m$}
\]
or
\[
\mathcal L_{\mS} = H \cup   \left(  \{(\beta+1) m_n, (\beta + 1) m_n- e\} + \mS \right), \;\; \hbox{\rm when $b=M$}.
\]

\end{enumerate}
\begin{enumerate}
\item [\rm{(II)}] Suppose that $b\not\in \{m,M\}$. Then 
\[
\mathcal L_{\mS} = H \cup   \left( \frac{e}{d}b + \mS \right).
\]
 \end{enumerate}
\end{thm}
\proof 
Let us prove \rm{(I)}. We first assume that $b = m$ and define $m_i'= m_i -b$ \textcolor{magenta}{for } $i\in \{1, \ldots, n\}$ and $m_i'' = \frac{m'_i}{d}$ with $d=\gcd(m'_1, \ldots, m'_n)$. Now, by Remark \ref{equivalent-ideal}.(1) and \ref{equivalent-ideal}.(3) we know that $I_{\tilde{\mS} } = I_{\mathcal T}$, where 
$\mathcal T = \langle (0,1), (m_1'', 1) , \ldots, (m_n'', 1) \rangle \subseteq \n^2.$
Since $m''_1<\ldots < m_n''$ is an arithmetic sequence of relative primes, we can apply Remark \ref{rem:resultado-Eva} to obtain a set of generators of the ideal $I_{\mathcal T} = I_{\tilde{\mS}} = \langle g_1, \ldots, g_s\rangle$ and then, Proposition \ref{prop-main-dS} to obtain $\mL_{\mS}$. In fact, if we set $l \equiv \frac{b-m_n+d}{d}\mod (n-1)$ with $l \in \{1,\ldots,n-1\}$, then
\begin{eqnarray*}
\mL_{\mS}  & = & \bigcup_{i=1}^s\, (\deg_{\mS}(g_i)+\mS) =  H \cup (\{\beta m_1 + m_i \mid 1\leq i \leq l \}+ \mS).
\end{eqnarray*}

Moreover, observe that for $i\geq 3$, then $m_1 + m_i = m_2 + m_{i-1}$. Thus, 
\[\beta m_1 + m_i = (\beta-1) m_1 + m_2 + m_{i-1} \in H.\] 
With this observation the above formula for $\mathcal L_{\mS}$ can be simplified as follows:
\begin{itemize}
\item If $l=1$ (or, equivalently, $d(n-1)$ divides $m_n - b$), then,
\[\mL_{\mS} = H \cup \left( (\beta+1) m_1 + \mS\right).\]
\item If $l\neq 1$, then, $\mL_{\mS} = H \cup \left( \{ (\beta+1) m_1, (\beta+1) m_1 + e\} + \mS\right).$
\end{itemize}
When $b = M$, we apply Remark \ref{equivalent-ideal}.(2) and the proof is analogue to that of $b = m$. 
\medskip

Now, let us prove \rm{(II)}. By Remark \ref{equivalent-ideal}.(1) we know that $I_{\tilde{\mS} } = I_{\mS_1}$ where
\[
\mathcal S_1 = \langle (0,1), (B, 1) , \ldots, ((n-1)B, 1), (c, 1)\rangle,
\]
being $c = \frac{b-m_1}{d}$ and $B:=\frac{e}{d}=\gcd (B, 2B, \ldots, (n-1)B)$.

Let us find explicit $\alpha_i \in \{1, \ldots, n\}$ such that
$B \cdot c = \sum_{i=1}^{n-1}\alpha_i \cdot i \cdot B$  with  $\sum_{i=1}^{n-1}\alpha_i  \leq B.$
We take $s \in \{1,\ldots,n-1\}$ such that $m_s < b < m_{s+1}$; then  $(s-1)B < c < sB$. Performing euclidean division we get
$c = \mu s + r$ with $1 \leq \mu < B$ and $r \in \{0,\ldots,s-1\}$. Then, $B c = \mu (sB) + (rB)$ and $\mu + 1 \leq B$.

By Remark \ref{rm:quitaruno} we have  
$I_{\mathcal S_1} = I_{\mathcal S_2} \cdot \mathbb K[x_1,\ldots,x_{n+1}] + \langle x_{n+1}^{B} - x_1^{B-\mu - 1} x_{r+1} x_{s+1}^{\mu}  \rangle,$
with $\mS_2 = \langle (0,1), (e, 1) , \ldots, ((n-1) e, 1) \rangle$. Moreover, applying Remark \ref{equivalent-ideal}.(3) we get $I_{\mS_2} = I_{\mS_3},$
with $\mS_3 = \langle (0,1), (1, 1) , \ldots,  (n-1, 1) \rangle \subseteq \n^2.$ 
Since $I_{\mS_3} = \langle x_i x_j - x_{i-1}x_{j+1} \mid 2\leq i \leq j \leq n-1\rangle,$
we can finally apply Proposition \ref{prop-main-dS} to obtain $\mathcal L_{\mS}$ from the set of generators of the ideal $I_{\tilde{\mS} }$. Thus,
$\mathcal L_{\mS} = \left( B\cdot b + \mS\right)  \cup \left( \{m_1 + \lambda e  \mid 2\leq \lambda \leq 2n-4 \} + \mS\right).$
\eproof

\medskip

We finish this section with an example illustrating Theorem \ref{prop:almostarithmetic-A}.

\begin{ex}
Let $\mS = \langle b, m_1, m_2, m_3, m_4, m_5\rangle$ be the numerical semigroup generated by  $b =7$, $m_1=17$, $m_2=20$, $m_3=23$, $m_4=26$ and $m_5=29$.
Note that $m_1 < \cdots < m_5$ is an arithmetic sequence of $n = 5$ relative primes, being $e=3$ the difference between two consecutive terms. We observe that $b \leq m_i$ for all $i \in \{1,\ldots,5\}$ and define
\[\begin{array}{ccc}
d=\gcd(m_1-b,e) = 1 & \hbox{ and } &  \beta = \left\lfloor \frac{m_n-b-d}{d(n-1)} \right\rfloor = 5,
\end{array}\]
and remark that $d(n-1)$ does not divide $m_n - b$. Then, by Theorem \ref{prop:almostarithmetic-A} we have 
\begin{eqnarray*}
\mL_{\mS} & = & (\{ 40,43,46,49,52\} + \mS) \cup (\{102,105 \} + \mS)\\
 &  = & \{ 40,43,46,49,52,102,105\} + \mS.
 \end{eqnarray*}
\end{ex}

\begin{cor} Let $\mS $ be the numerical semigroup generated by the almost arithmetic sequence $\mathcal A = \{ m_1, \ldots, m_n, b\} \subseteq \n$ as in \eqref{eq-almostarithmetic}. Put
 $M:=\max \mA,\, m:=\min \mA$ and $d := \gcd(e,b-m_1)$. Then
\begin{enumerate}
\item For $b \in \{m,M\}$ we get
\[{\rm c_{\rm eq}}(\mS) = \left\lceil \frac{M-m -d-1}{d(n-1)} \right\rceil.\]
\item For $m_1 < b < m_n$ we get
\[{\rm c_{\rm eq}}(\mS) = \frac{e}{d}.\]
\end{enumerate}
\end{cor}

\proof
Following the lines of the proof of Theorem \ref{prop:almostarithmetic-A} one observes that the maximum degree in a minimal set of generators of $I_{\tilde{S}}$ is $\frac{e}{d}$ if $m_1 < b < m_n$, or $\left\lfloor \frac{M-m-d}{d(n-1)} \right\rfloor + 1$  when $b \in \{m,M\}$. The result follows from Theorem \ref{th:catenary}. 
\eproof

\section{When is $\mL_\mS$ a principal ideal?}
\label{sec:principal}

Whenever $\mS = \langle \ba_1,\ldots,\ba_n\rangle \subseteq \zz^m \oplus T$ is a reduced monoid such that $\mL_\mS = e + \mS$ for some $e \in \mS$, we have that $x \in \mS$ if and only if $e + x \in \mL_\mS$. When $\mS$ is a numerical semigroup, the previous trivial observation implies that, in particular, the maximum element not in $\mL_\mS$ and $F(\mS)$, the Frobenius number of $\mS$, are closely related. Indeed, $\max \{b \in \zz \,\vert \,  b \notin \mL_\mS\} = e + F(\mS)$. 
This is one of the reasons why  it could be interesting to characterize numerical semigroups such that $\mL_\mS$ is a principal ideal.

When $\mS = \langle a_1,a_2,a_3\rangle$ is a three-generated numerical semigroup, then
$\mL_\mS$ is a principal ideal (see \cite{chapman:2011}). In  Corollary \ref{cor:numericos}, we provided another proof of the same fact. The idea in our proof is that 
$I_{\tilde{\mS}}$ is a height one ideal and, thus, it is principal. As a consequence, this proof can be generalized to reduced monoids $\mS = \langle \ba_1,\ldots,\ba_n\rangle \subseteq \zz^m  \oplus T$ as far as $I_{\tilde{\mS}}$ is a height one ideal (see also Proposition \ref{prop-main-dS}). 
However, this is not the only situation in which $\mL_\mS$ is a principal ideal. In Corollary \ref{cor:uniqueBetti} we provide a family of numerical semigroups such that $\mL_\mS$ is a principal ideal.  This family includes the one of three-generated numerical semigroups.

We begin with a proposition which follows from Proposition \ref{prop-main-dS}:
\begin{prop}
\label{prop:lsprincipal}
Let $\mS \subseteq \zz^m \oplus T$ be a finitely generated reduced monoid and take $\{g_1, \ldots g_r\}$ a binomial generating set of $I_{\tilde{\mS}}$. Then,
$\mL_{\mS}$ is a principal ideal if and only if there exists  $i \in \{1,\ldots,r\}$ such that $\deg_\mS(g_j) \in \deg_\mS(g_i) + \mS$ for all $j \in \{1,\ldots,r\}$.
\end{prop}

We observe that the above  condition on $\mS$-degrees can be restated as fo\-llows: if one considers $\leq_\mS$ the partial order $y \leq_\mS z$ if and only if $z - y \in \mS$, then the set of $\mS$-degrees of the generators of $I_{\tilde{\mS}}$ has a minimum element. 
This condition for $\tilde{\mS}$ is slightly more general than the one of being an {\it affine monoid with one Betti minimal element}, explored in  \cite{garcia-sanchez:2019}.  In this section, we build on  some ideas of \cite[Section 7]{garcia-sanchez:2019}.

Now we describe $\mL_\mS$ for a particular family of numerical semigroups.

\begin{prop}\label{prop:descuniquebetti} Let $\mS = \langle b, b + t m_1, \ldots, b + t m_n \rangle$ be a numerical semigroup, where $b, t \in \zz^+$, $n \geq 2$ and  $m_i = f_i \prod_{j \in \{1,\ldots,n\} \atop j \neq i} c_j$; being
\begin{itemize}
\item[{\rm (a)}]  $c_1,\ldots,c_n \in \n$ pairwise relatively prime,
\item[{\rm (b)}]  $\gcd(f_i,c_i) = 1$ for all $i \in \{1,\ldots,n\}$, 
\item[{\rm (c)}]  $m_n > m_i$ for all $i \in \{1,\ldots,n-1\}$, and
\item[{\rm (d)}]   $f_n = 1$.
\end{itemize}
Then $\mL_\mS = \bigcup_{i = 1}^{n-1} \left(c_i (b + tm_i) + \mS\right)$.
\end{prop}
\proof We will make use of Proposition \ref{prop-main-dS}. For this purpose, we are obtaining a generating set for $I_{\tilde{\mS}}$. By Remark \ref{equivalent-ideal}.(1) we have  $I_{\tilde{\mS}} = I_T$, where $T = \langle (0,1),(m_1,1),\ldots,(m_n,1) \rangle$.  We observe that $\gcd(m_1,\ldots,m_n) = 1$, and for all $i \in \{1,\ldots,n-1\}$ we get 
\[ \gcd(m_1,\ldots,m_{i-1},m_{i+1},\ldots,m_n) m_i = c_i m_i = f_i c_n m_n, \] and $f_i c_n < c_i$ (because $m_i < m_n$). Thus, applying Remark \ref{rm:quitaruno}, we have  $I_T = \langle x_{i+1}^{c_i} - x_1^{c_i - f_i c_n } x_n^{f_i c_n} \, \vert \, 1 \leq i \leq n-1 \rangle$. Since ${\rm deg}_{\mS}(x_{i+1}^{c_i}) = c_i (b + tm_i)$ for $i \in \{1,\ldots,n-1\}$, by Proposition \ref{prop-main-dS} we are done. \eproof 

In the proof of Proposition \ref{prop:descuniquebetti} we obtain a minimal set of generators of $I_{\tilde{\mS}}$. Hence, applying Theorem \ref{th:catenary} we get:
\begin{cor}
Let $\mS = \langle b, b + t m_1, \ldots, b + t m_n \rangle$ be a numerical semigroup, where $b, t \in \zz^+$, $n \geq 2$ and  $m_i = f_i \prod_{j \in \{1,\ldots,n\} \atop j \neq i} c_j$; being
\begin{itemize}
\item[{\rm (a)}]  $c_1,\ldots,c_n \in \n$ pairwise relatively prime,
\item[{\rm (b)}]  $\gcd(f_i,c_i) = 1$ for all $i \in \{1,\ldots,n\}$, 
\item[{\rm (c)}]  $m_n > m_i$ for all $i \in \{1,\ldots,n-1\}$, and
\item[{\rm (d)}]   $f_n = 1$.
\end{itemize}
Then, ${\rm c_{\rm eq}}(\mS) = \max \{ c_i \mid 1 \leq i \leq n-1\}.$
\end{cor}

Now, we apply Proposition \ref{prop:descuniquebetti} to the subfamily of the numerical semigroups, which corresponds to setting $f_i = 1$ for all $i \in \{1,\ldots,n\}$. Hence, the semigroup  $\mS$ belongs to the so-called {\it shifted family} of $\mS' = \langle m_1,\ldots, m_n \rangle$, where $\mS'$ is a {\it numerical semigroup with a unique Betti element}; we refer the reader to \cite{garcia-sanchez:2012} for more on semigroups with a unique Betti element.

\begin{cor}\label{cor:uniqueBetti} Let $\mS = \langle b, b + t m_1, \ldots, b + t m_n \rangle$ be a numerical semigroup, where $b, t \in \zz^+$ and  $m_i = \prod_{j \in \{1,\ldots,n\} \atop j \neq i} c_j$; being  $c_1 > \cdots > c_n \geq 2$ pairwise relatively prime integers.
Then, $\mL_\mS = c_{n-1} (b + tm_{n-1}) + \mS$.
\end{cor}
\proof Clearly the hypotheses of Proposition \ref{prop:descuniquebetti} are satisfied with $f_i = 1$ for all $i \in \{1,\ldots,n\}$. Set $D_i := c_i (b + t m_i)$ for all $i \in \{1,\ldots,n-1\}$. To conclude, it suffices to prove that $D_i \in D_{n-1}  + \mS$ or, equivalently, that $D_i - D_{n-1} \in \mS$ for all $i \in \{1,\ldots,n-2\}$. Take $i \in \{1,\ldots,n-2\}$, we have that  
\[ \begin{array}{lll} D_{i} - D_{n-1} & = & (c_i - c_{n-1})b + t (c_i m_i - c_{n-1}m_{n-1}) =  (c_i - c_{n-1})b \in \mS \end{array}. \]
\eproof 

Let us illustrate Corollary \ref{cor:uniqueBetti}  with an example.

\begin{ex} Let $\mS = \langle 17, 29, 37, 47 \rangle$, which satisfies the hypotheses of Proposition \ref{prop:lsprincipal}, with $b = 17,\, t = 2,\, n = 3,\, c_1 = 5,\, c_2 = 3$ and $c_3 = 2$. Thus, $\mL_\mS = (3 \cdot 37) + \mS = 111 + \mS$. 
Indeed, as we proved in Proposition \ref{prop:descuniquebetti} and Corollary \ref{cor:uniqueBetti}, $I_{\tilde{\mS}} = \langle g_1,g_2 \rangle$ with $g_1 = x_3^3-x_1 x_4^2$ and $g_2 = x_2^5 - x_1^3x_4^2$ and we have  $\deg_\mS(g_2) \in \deg_\mS(g_1) + \mS$, because $\deg_{\mS}(g_1) = 3 \cdot 37 = 111,\ \deg_{\mS}(g_2) = 5 \cdot 29 = 145= 111 + 2 \cdot 17 \in 111 + \mS$. 
Moreover, since the Frobenius number of $\mS$ is $F(\mS) = 107$, we get $\max\{b \in \zz \, \vert \, b \notin \mL_\mS\} = 111 + 107 = 218$.  
\end{ex}

One could build further families of numerical semigroups such that $\mL_\mS$ is a principal ideal by  choosing appropriate values of $f_1,\ldots,f_{n-1}$ in Proposition~\ref{prop:descuniquebetti}.

We observe that  Corollary  \ref{cor:uniqueBetti} includes the case of three generated numerical semigroups and, hence, generalizes the formula obtained in Corollary \ref{cor:numericos}. Indeed, the numerical semigroup $\mS = \langle a_1, a_2, a_3 \rangle$ with $a_1<a_2<a_3$ corres\-ponds to $b = a_1$, $n = 2$, $t = \gcd(a_2-a_1,a_3-a_1)$, $m_1 = c_2 = (a_2-a_1)/t$ and $m_2 = c_1 = (a_3 - a_1)/t$ and, in this context, we have   \[ \mL_\mS = c_{1} (b + tm_{1}) + \mS = (a_2 (a_3-a_1) / \gcd(a_2-a_1,a_3-a_1)) + \mS.\]


\section{Computational considerations} 
\label{sec:NP}

Let $\mS = \langle a_1,\ldots,a_n \rangle$ be a numerical semigroup, as we saw in Corollary \ref{cor:numericos}, then $\mL_\mS = \emptyset$ if and only if $n \leq 2$. Thus, when $n \geq 3$, by Corollary \ref{cor:chapman} if follows that $\n \setminus \mL_\mS$ is a finite set. Hence, for $n \geq 3$ the integer $F_{2,\ell} = \max\{b  \in \zz \mid b \notin \mL_{\mS} \}$ is well defined. 

 The goal of this short section is to show that the problem of computing the largest element in $\zz \setminus \mL_\mS$ is an $\mathcal{NP}$-hard problem, under Turing reductions. 
 
 In \cite{alfonsin:1996} (see also \cite[Theorem 1.3.1]{alfonsin:2005}), Ram\'irez Alfons\'in proves that the problem of determining the Frobenius problem is $\mathcal{NP}$-hard. His proof consists of a Turing reduction from the Integer Knapsack Problem (IKP), which is well-known to be an $\mathcal{NP}$-complete problem (see, e.g., \cite[page 376]{PScomplexity}). The ${\rm IKP}$ is a decision problem that receives as input $(a_1,\ldots,a_n) \in \n^n$, $t  \in \n$ and asks if there exist $x_1,\ldots,x_n \in \n$ such that $\sum_{i = 1}^n x_i a_i = t$. We define here a related decision problem, we call this problem ${\rm IKP}_{2,\ell}$: 
 \begin{itemize} 
 \item {\bf Input:} $(a_1,\ldots,a_n) \in \n^n,\, t \in \n$, and 
 \item {\bf Question:} do there exist distinct $(x_1,\ldots,x_n), (y_1,\ldots,y_n) \in \n^n$ such that $\sum_{i = 1}^n x_i a_i = \sum_{i = 1}^n y_i a_i = t$ and $\sum_{i = 1}^n x_i = \sum_{i = 1}^n y_i$?
\end{itemize}
Observe that 
\begin{center}
${\rm IKP}((a_1,\ldots,a_n),t) =$ {\sc True} $\Leftrightarrow\ {\rm IKP}_{2,\ell}((a_1,\ldots,a_n,a_1,\ldots,a_n),t) =$~{\sc True}, \end{center}
implies that ${\rm IKP}_{2,\ell}$ is an $\mathcal{NP}$-hard problem.

Moreover, a careful inspection of the proof of \cite[Theorem 1.3.1]{alfonsin:2005} shows that if we replace ${\rm IKP}$ by ${\rm IKP}_{2,\ell},$ and $F(\langle a_1,\ldots,a_n \rangle)$ by $F_{2,\ell}(\langle a_1,\ldots,a_n \rangle)$ the proof also holds. This fact together with the $\mathcal{NP}$-hardness of ${\rm IKP}_{2,\ell}$ yields the following:
\begin{prop}\label{pr:nphard}Let $\mS = \langle a_1,\ldots,a_n \rangle$ be a numerical semigroup with $n \geq 3$. The problem of computing $F_{2,\ell}(\mS) = \max\{b \in \zz \mid b \notin \mL_{\mS}\}$ is $\mathcal{NP}$-hard.
\end{prop}

We finally remark that one can define 
\[ \begin{array}{cll} F_i(\mS) & = & {\rm max}\{b \in \zz \, \vert \, b \text{  has  not  } i \text{\ factorizations}\}, {\text \ and} \\
F_{i,\ell}(\mS) & = & {\rm max}\{b \in \zz \, \vert \, b \text{  has not } i \text{ factorizations of the same length}\} \end{array} \]  and, following the same argument presented here, one can prove that the computational problem of computing $F_i(\mS)$  or $F_{i,\ell}(\mS)$ for bounded values of $i$ are all $\mathcal{NP}$-hard.

\section*{Acknowledgements}
We wish to thank M.A. Moreno-Fr\'ias, who presented the results of \cite{GarciaGarcia:2019} in
the seminar GASIULL at Universidad de La Laguna and introduced and
encouraged us to work on this topic.

We also want to thank the anonymous referee for his/her insightful comments. In particular, the referee suggested to relate the results in an earlier version of the paper with the equal catenary degree. This suggestion gave rise to Section~\ref{sec:catenaryDegree}. Moreover, in this earlier version, we considered the problem of factorizations in affine monoids and the referee suggested to tackle the same problem in the (slightly) more general context of abelian, cancellative and finitely generated monoids. Following this suggestion we could obtain more general results in Sections~\ref{Sec:Apery} and \ref{sec:samelenght}.

\section*{Fundings}

This work was partially supported by the Spanish MICINN PID2019-105896GB-I00 and MASCA (ULL Research Project).

\section*{Conflict of interest}
On behalf of all authors, the corresponding author states that there is no conflict of interest.

\bibliographystyle{plain}
\bibliography{ArxivVersion}

\end{document}